\documentclass[english]{amsart}
\usepackage{amsopn}                    
\usepackage{amsmath}
\usepackage{amssymb}
\usepackage{stmaryrd}
\usepackage[all,2cell,emtex]{xy}
\usepackage[mathscr]{euscript}
\usepackage{mathrsfs}
\usepackage{bbm}
\DeclareFontFamily{OT1}{pzc}{}
\DeclareFontShape{OT1}{pzc}{m}{it}%
             {<-> s * [1,150] pzcmi7t}{}
\DeclareMathAlphabet{\mathpzc}{OT1}{pzc}%
                                 {m}{it}

\theoremstyle{plain} 
\newtheorem{thm}{Theorem}[section]
\newtheorem{prop}[thm]{Proposition}
\newtheorem{lemma}[thm]{Lemma}
\newtheorem{cor}[thm]{Corollary}

\theoremstyle{remark}
\newtheorem{rem}[thm]{Remark}

\theoremstyle{definition}
\newtheorem{definition}[thm]{Definition}
\newtheorem{paragr}[thm]{}

\theoremstyle{plain} 
\swapnumbers

\numberwithin{equation}{thm}

\renewcommand{\mathcal}{\mathpzc}
\renewcommand{\mathbb}{\mathbbm}

\renewcommand{\leq}{\leqslant}
\renewcommand{\geq}{\geqslant}

\newcommand{\pref}[1]{{\widehat{ #1 }}}

\newcommand{\To}{\longrightarrow}
\newcommand{\cats}{\Delta}

\newcommand{\Hom}{\operatorname{\mathrm{Hom}}}
\newcommand{\sHom}{\operatorname{\mathpzc{Hom}}}

\newcommand{\nerf}{N}

\newcommand{\op}[1]{{#1}^{\mathit{op}}}

\newcommand{\real}[1]{|#1|_J}
\newcommand{\sing}[1]{\mathit{Sing}_J(#1)}

\newcommand{\ho}{\operatorname{\mathbf{Ho}}}

\newcommand{\C}{\mathcal{C}}

\newcommand{\smp}[1]{ \Delta[#1]}

\newcommand{\derR}{\mathbf{R}}

\newcommand{\bord}{\partial}

\def\TO#1{\mathrel{\hbox to #1pt{\rightarrocharReedyfibrantsegaloperwfill}}}
\def\OT#1{\mathrel{\hbox to #1pt{\leftarrowfill}}}

\def\limproj{\mathop{\oalign{\rm lim\cr
\hidewidth$\longleftarrow$\hidewidth\cr}}}%
\def\limind{\mathop{\oalign{\rm lim\cr
\hidewidth$\longrightarrow$\hidewidth\cr}}}%
\renewcommand{\varinjlim}{\limind}%
\renewcommand{\varprojlim}{\limproj}%

\newcommand{\Wpr}{\mathsf{W}}

\renewcommand{\to}{\To}
\newcommand{\todouble}{\xymatrixcolsep{1pc}\xymatrix{\ar@<.5ex>[r]\ar@<-.5ex>[r]&}}
\newcommand{\todoubleop}{\xymatrixcolsep{1pc}\xymatrix{\ar@<.5ex>[r]&\ar@<.5ex>[l]}}

\renewcommand{\hookrightarrow}{{\hskip -1.5pt\raise 1.5pt\vbox{\xymatrixcolsep{.9pc}\xymatrix{\ar@{^{(}->}[r]&}}\hskip -3.5pt}}

\newcommand{\intcoin}[5]{\raise 12pt\vbox{\xymatrixcolsep{.9pc}\xymatrixrowsep{.7pc}\xymatrix{%
\scriptstyle #1\ar[r]^{\scriptscriptstyle #5}\ar[d]_{\scriptscriptstyle #4}&\scriptstyle #3\\\scriptstyle #2}}}

\renewcommand{\xrightarrow}[1]{{\hskip -2.5pt\xymatrixcolsep{1.7pc}\xymatrix{\ar[r]^{#1}&}\hskip -2.5pt}}

\newcommand{\set}{\mathpzc{Set}}
\newcommand{\sset}{\mathpzc{sSet}}
\newcommand{\dset}{\mathpzc{dSet}}
\newcommand{\sdset}{\mathpzc{sdSet}}

\newcommand{\oper}{\mathpzc{Operad}}
\newcommand{\preoper}{\mathpzc{PreOper}}
\newcommand{\segoper}{\mathpzc{SegOper}}
\newcommand{\map}{\mathpzc{hom}}

\newcommand{\spine}{\mathrm{Sc}}
\newcommand{\bext}{{\partial^\mathit{ext}\Omega}}

\title{Dendroidal Segal spaces and $\infty$-operads}

\date{\today}

\author[D.-C. Cisinski]{Denis-Charles Cisinski}
\address{Universit\'e Paul Sabatier\\
Institut de Math\'ematiques de Toulouse\\
118~route de Narbonne\\
31062~Toulouse cedex~9\\France}
\email{denis-charles.cisinski@math.univ-toulouse.fr}
\urladdr{http://www.math.univ-toulouse.fr/~dcisinsk/}

\author[I. Moerdijk]{Ieke Moerdijk}
\address{Radboud Universiteit Nijmegen\\
Institute for Mathematics, Astrophysics, and Particle Physics\\
Heyendaalseweg~135, 6525~AJ~Nijmegen\\
The~Netherlands}
\email{i.moerdijk@math.ru.nl}
\urladdr{http://www.math.uu.nl/people/moerdijk/}


\begin{document}
\begin{abstract}
 We introduce the dendroidal analogues of the notions of complete Segal space and
of Segal category, and construct two appropriate model categories for which each of these notions
corresponds to the property of being fibrant. We prove that these two model categories are
Quillen equivalent to each other, and to the monoidal model category for $\infty$-operads which we
constructed in an earlier paper. By slicing over the monoidal unit objects in these model
categories, we derive as immediate corollaries the known comparison results between Joyal's
quasi-categories, Rezk's complete Segal spaces, and Segal categories.
\end{abstract}
\maketitle
\tableofcontents

\section*{Introduction}
The category of dendroidal sets is an extension of that of simplicial sets, suitable for constructing nerves,
not just of categories but also of (coloured) operads. It was introduced with this purpose, and with the aim
of giving an inductive definition of weak higher categories, in \cite{dend1,dend2}.
This category $\dset$ of dendroidal
sets carries a symmetric monoidal closed structure which is closely related to the Boardman-Vogt
tensor product of operads, and the inclusion of the category $\sset$ of simplicial sets into $\dset$ can in
fact be identified with the forgetful functor, from the slice (or comma) category of dendroidal sets over
the unit $\eta$ of the monoidal structure  back to dendroidal sets, via an explicit isomorphism of categories
\begin{equation}\label{intro:compdsetoveretasset}
\dset/\eta=\sset
\end{equation}

Dendroidal sets carry  a very rich homotopical structure, which we began to explore in \cite{dend3}.
For example, there is a monoidal Quillen model structure on $\dset$, whose fibrant objects include
all nerves of operads. In fact, these fibrant objects can be thought of as simple combinatorial models
of the notion of operad-up-to-homotopy or ``$\infty$-operad''. Like any Quillen model structure,
this model structure on dendroidal sets induces another model structure on any slice category.
Under the identification $\dset/\eta=\sset$  this induced model structure can be shown to coincide with
the Joyal model structure on simplicial sets, whose fibrant objects are most commonly known under the
name ``$\infty$-categories'' (and are also referred to as quasi-categories,  weak Kan complexes, or inner
Kan complexes \cite{joyal,lurie,BV}).  

These $\infty$-categories model a notion of category-up-to-homotopy.
Other ways of modelling such a notion have occurred in the literature,
including the theory of Segal categories~\cite{simpson,bergner} and of
complete Segal spaces~\cite{Rezk}. The latter two concepts are both based on the
much older observation that a simplicial set X is the nerve of a category if and only if the canonical map 
\begin{equation}\label{intro:Segalmap}
X_n \To  X_1 \times_{X_0}\dots \times_{X_0} X_1
\end{equation}
sending a simplex to its one-dimensional ribbons, is an isomorphism. Indeed,
Simpson and Rezk both base their theories on bisimplicial sets X for which the map \eqref{intro:Segalmap}
is a weak equivalence of simplicial sets (and replacing the fibred product on the right
hand side by its homotopy version). Building on the work of Simpson and Rezk,
the relation between these different ways of modelling categories-up-to-homotopy was
recently made precise through the work of  Bergner, Joyal and Tierney, and Lurie.
Indeed, Simpson's Segal categories, Rezk's complete Segal spaces, and Joyal's $\infty$-categories
all arise as the fibrant objects in a specific Quillen model category structure, and these different
model category structures have now been related to each other by explicit Quillen
equivalences \cite{bergner,lurie}.
Moreover, they are all Quillen equivalent to the model category of simplicial categories discovered by
Bergner~\cite{bergner},
thus providing a strictification or rigidification result for each of these notions of category-up-to-homotopy. 

The goal of this paper and its sequel \cite{dend5}
is to develop  analogous theories of Segal operads (rather than categories) and complete
dendroidal (rather than simplicial) Segal spaces, to relate these to each other and to dendroidal sets
via Quillen equivalences, and to prove a strictification result for each of them by relating them to
simplicial operads. By a simple slicing procedure like in \eqref{intro:compdsetoveretasset},
the earlier results just mentioned for
categories-up-to-homotopy can all be recovered from our results, which can in this sense be said to be
more general.

In more detail, then, we will consider the category $\sdset$ of simplicial objects in dendroidal sets,
or what is the same, dendroidal spaces. We will define a Segal type condition on the objects of this
category, based on an extension to trees of ``the union of $1$-dimensional ribbons in an $n$-simplex'' to
which we will refer as the Segal core of a tree. In Section~\ref{section:dendsegalspaces},
we will establish a closed model category
structure on $\sdset$ whose fibrant objects satisfy a tree-like Segal condition involving these Segal
cores, and a completeness condition like the one of Rezk, and prove (Corollary~\ref{inclQuillenequiv})
that this model category is Quillen equivalent to our earlier model category structure on dendroidal
sets~\cite{dend3}. The definitions and proofs of these results are based on  some elementary observations
about these Segal cores presented in Section~\ref{section:cores},
and on a characterization of weak equivalences between
$\infty$-operads as maps which are ``essentially surjective and fully faithful'' in a suitable sense
(Theorem~\ref{charactwoperequiv1}).
The proof also exploits the hybrid nature of the objects of $\sdset$, which can be viewed
alternatively as simplicial objects in one category or as dendroidal objects in another. In fact, the first
view point is taken in Section \ref{section:locconstant},
while the second viewpoint underlies the notion of complete
dendroidal Segal space and the formulation of the main equivalence \ref{inclQuillenequiv}.
The relation between these two view points is most clearly expressed by Theorem~\ref{comploccstRezk}
which equates two seemingly different model category structures. 

Again using the Segal cores, we define the notion of a Segal operad in Section~\ref{segaloperads}.
These Segal operads will then be shown to be the fibrant objects for a model category structure on a full 
subcategory of the category $\sdset$ of dendroidal spaces,  the category of so-called Segal pre-operads
(Theorems \ref{thm:modcatSegaloperads} and \ref{charReedyfibrantsegaloper}).
Using most if not all of the earlier results, we will then be able to show that this model category with
Segal operads as fibrant objects is Quillen equivalent to the model category having complete dendroidal
Segal spaces as fibrant objects (Theorem~\ref{QuillenequivSegal}),
and hence also Quillen equivalent to the original model category of dendroidal sets.

We believe these results are of interest in themselves, and because they generalize important
classical results from the simplicial-categorical context to the dendroidal-operadic one. In addition,
they will all be used in our proof of the strictification theorem for $\infty$-operads presented in \cite{dend5}.

\section{Preliminaries}

We begin by recalling the basic definitions related to dendroidal sets; see \cite{dend1,dend2,dend3}.
The starting point is a category $\Omega$ of trees. Its objects are finite (non-planar) trees.
These trees have internal edges (between vertices) and external ones (attached to just one vertex);
the root is one such external edge, the others are called ``leaves'', or ``input edges''.
Each such tree freely generates a coloured operad, and the arrows in $\Omega$ are the maps
between these operads. Thus, by definition, $\Omega$ is a full subcategory of the category
of (symmetric coloured) operads.

Each natural number $n\geq 0$ defines a linear tree with $n$ vertices and $n+1$ edges, the input edge
being labelled $0$, and the output or root edge labelled $n$. The corresponding coloured operad
is the category defined by the linear order $0\leq \cdots \leq n$. Thus the simplicial category $\Delta$
is a full subcategory of $\Omega$, and we denote the inclusion by
$$i:\Delta\To \Omega\ ,  \quad [n]\longmapsto n=i[n]\, .$$

The category $\dset$ of \emph{dendroidal sets} is by definition the category of presheaves (i.e.
contravariant $\set$-valued functors) on $\Omega$, just like the category of \emph{simplicial sets}
is that of presheaves on $\Delta$. The inclusion functor $i$ induces a pair of adjoint functors
$$i_!: \sset\rightleftarrows\dset:i^*$$
where $i^*$ is the restriction along $i$ and $i_!$ is its fully faithful left adjoint ($i^*$
also has a fully faithful right adjoint $i_*$).

We will write $\Omega[T]$ for the dendroidal set represented by a tree $T$. With the similar notation
$\Delta[n]$ for representable simplicial sets, we thus have
$$i_!\Delta[n]=\Omega[n]\, ,$$
and this identification determines $i_!$ uniquely up to unique isomorphism (as colimit
preserving functor).

There is a natural identification of $\Delta$ with the slice category $\Omega/i[0]$, and this
leads to an identification
$$\sset=\dset/\eta$$
where $\eta=\Omega[0]$. Under this identification, the functor $i_!$ corresponds to the forgetful
functor $\dset/\eta\to \dset$.

The full embedding of $\Omega$ into (coloured) operads gives an adjoint pair
$$\tau_d:\dset\rightleftarrows\oper:\nerf_d\, ,$$
where the right adjoint $\nerf_d$ is called the \emph{dendroidal nerve}.
These functors restrict to the usual nerve of a small category and its left adjoint.

The category of dendroidal sets carries a symmetric monoidal closed structure, denoted
by $\otimes$ and $\sHom$. Its unit object is the representable dendroidal set $\eta=\Omega[0]$.
This structure is compatible with the product of simplicial sets as well as with
the Boardman-Vogt tensor product
of operads, in the sense that, for any simplicial sets $M$ and $N$, and for any dendroidal sets $X$ and $Y$,
we have natural identifications
$$i_!(M\times N)=i_!(M)\otimes i_!(N)\quad\text{and}
\quad \tau_d(X\otimes Y)=\tau_d(X)\otimes_{\mathit{BV}} \tau_d(Y)\, .$$

We now recall some of the main combinatorial properties in the category of dendroidal sets;
see \cite{dend1,dend2,dend3} for more details.

Just like for the simplicial category $\Delta$, the arrows in $\Omega$
are generated by elementary arrows. These are faces and degeneracies like for $\Delta$, together
with the isomorphisms (the only isomorphisms in $\Delta$ are the identities).
In particular, for a tree $T$, we may define $\bord \Omega[T]$ as the maximal proper subobject
of $\Omega[T]$, or, equivalently, as the union of the all the images of the elementary
face maps $\Omega[S]\To \Omega[T]$. We refer to $\bord\Omega[T]$
as the \emph{boundary} of $\Omega[T]$. The saturation of the set of boundary inclusions
$\bord \Omega[T]\to \Omega[T]$ (i.e. the closure under transfinite composition, pushout, and retract)
gives rise to the class of \emph{normal monomorphisms}. The normal monomorphisms
can also be characterized as the monomorphisms of dendroidal sets $u:X\To Y$
such that, for any tree $T$ in $\Omega$, the action of $\mathrm{Aut}(T)$ on the set $Y(T)-u(X(T))$
is free. A dendroidal set $X$ is \emph{normal} if the map $\varnothing\To X$ is a normal
monomorphism. We will often use the following property: given any morphism of
dendroidal sets $X\To Y$, if $Y$ is normal, then so is $X$.

For an internal edge $e$ in a tree $T$, we denote by $T/e$ the tree obtained from $T$ by contracting
the edge $e$. Then there is an elementary face map
$$\partial_e:T/e\To T$$
in $\Omega$. Face maps of this shape are called \emph{inner} or
\emph{internal}. We write $\Lambda^e[T]$
for the maximal subobject of $\Omega[T]$ which does not contain the image of the
internal face $\partial_e:\Omega[T/e]\To \Omega[T]$ (equivalently, $\Lambda^e[T]$
may be described as the union of all the images of the elementary faces $\Omega[S]\To \Omega[T]$
which do not factor through $\partial_e$).
We refer to $\Lambda^e[T]$ as the \emph{inner horn} of $\Omega[T]$ associated to $e$.
The class of \emph{inner anodyne extensions} is defined to be the closure of
the set of inner horn inclusions by transfinite composition, pushout, and retract.

A morphism of dendroidal sets $X\To Y$ is called an \emph{inner fibration} if it has the right
lifting property with repect to the set of inner horn inclusions $\Lambda^e[T]\To \Omega[T]$.
A dendroidal set $X$ is called an \emph{inner Kan complex}, or an \emph{$\infty$-operad}
if the map from $X$ to the terminal object of $\dset$ is an inner fibration.
Under the identification $\dset/\eta=\sset$, the $\infty$-operads which admit a (necessarily
unique) map to $\eta$ are precisely the $\infty$-categories (quasi-categories) of Joyal.
We may now formulate one of the main results of \cite{dend3}:

\begin{thm}\label{cmfdendsetjoyalike}
The category $\dset$ of dendroidal sets carries a cofibrantly generated model category structure,
whose cofibrations are the normal monomorphisms, and whose fibrant objets are the
$\infty$-operads. This structure is left proper and monoidal (i.e. compatible with tensor product).
The induced model category
structure on $\dset/\eta=\sset$ corresponds to the Joyal model structure on $\sset$.
\end{thm}

The weak equivalences of the model structure above are called the \emph{weak operadic
equivalences}. This model category structure on $\dset$ will be referred to as the
\emph{model category structure for $\infty$-operads}.


\section{Segal cores}\label{section:cores}

\begin{paragr}
We recall that, for each $n\geq 0$, the $n$th corolla $C_n$
is defined as the smallest rooted tree with one vertex and $n$ leaves.
\begin{equation}\label{corolla}
\begin{split}
\xymatrix@R=10pt@C=12pt{
&&&&&&\\
&*=0{}&&*=0{}\ar@{}[rrr]|{\cdots\cdots}&&&*=0{}&\\
C_n=&&&*=0{\bullet}\ar@{-}[u]_{a_2}\ar@{-}[ull]^{a_1}\ar@{-}[urrr]_{a_n}&&&\\
&&&*=0{}\ar@{-}[u]^a&&&
}
\end{split}\end{equation}
In general, we say that a face map $F\To T$ is a \emph {subtree}
if $F\To T$ is a composition of external faces.
In other words, a face map $F\To T$ is a subtree if $F$ is obtained by
successively pruning away top vertices, or pruning away root vertices
which have only one internal edge attached to them.
\begin{equation}
\begin{split}\xymatrix@R=5pt@C=6pt{
&&&&&&&&&\\
&&&&*=0{}\ar@{-}[rr]&\ar@{}[d]|F&*=0{}\ar@{-}[ll]&&&\\
&&&&&*=0{}\ar@{-}[ul]*=0{}\ar@{-}[ur] *=0{}\ar@{-}[d]&&&&\\
T=&&&&&&&&&\\
&&&&&*=0{}\ar@{-}[uuuullll]*=0{}\ar@{-}[uuuurrrr] *=0{}\ar@{-}[d]&&&&\\
&&&&&&&&&
}\end{split}
\end{equation}
\end{paragr}

\begin{definition}\label{def:spines}
Given a tree $T$ with at least one vertex, we define its \emph{Segal core} $\spine[T]$
as the subobject of $\Omega[T]$ defined as the union of all the images of those
maps $\Omega[C_n]\To \Omega[T]$ corresponding to subtrees of shape $C_n\To T$.
Remark that, up to isomorphism, such a map $C_n\To T$
is completely determined by the vertex of $T$ in its image, so we can write
$$\spine[T]=\bigcup_{v}\Omega[C_{n(v)}]\, ,$$
where the union is over all the vertices of $T$, and $n(v)$ is the number of input edges at $v$.

If $T=[0]$ is the tree with no vertices (so that $\Omega[T]=\Omega[0]=\eta$ is the
unit object of the Boardman-Vogt tensor product), it will be convenient to define $\spine[T]=\eta$.
\end{definition}

Recall from \cite[Paragraph 1.2]{dend3} that an (elementary) face $S\to T$ of a tree $T$
is called \emph{outer} or \emph{external} if $S$ is obtained from $T$
by pruning away an external vertex, i.e. a vertex with exactly one
inner edge attached to it.

\begin{definition}\label{def:externalboundary}
Let $T$ be a tree. The \emph{external boundary} of $\Omega[T]$ is the
subobject $\bext[T]$ of $\Omega[T]$ obtained as the union of all the
external faces of $T$.
\end{definition}

\begin{prop}\label{segalcoreanodyne}
For any tree $T$, the inclusion $\spine[T]\To\Omega[T]$
is inner anodyne.
\end{prop}

\begin{proof}
Note that, if $T$ has at most one vertex, then this inclusion is an isomorphism, while, if
$T$ has exactly two vertices, this is an inner horn.
So we may assume that $T$ has at least three vertices.

If $T$ has $N$ vertices, then $\Omega[T]$ has a natural filtration
by subobjects
$$\Omega[T]_1\subset\Omega[T]_2\subset \ldots \subset\Omega[T]_{N-1}\subset\Omega[T]_N=\Omega[T]\, ,$$
where, for $1\leq n\leq N$,
$$\Omega[T]_n=\bigcup_{\xi}\Omega[F_\xi]$$
is the union over all subtrees $F_\xi$ of $T$ with at most $n$ vertices.
Notice that, by definition, we have:
$$\Omega[T]_1=\spine[T]\quad \text{and} \quad
\Omega[T]_{N-1}=\bext[T]\, .$$
By virtue of \cite[Lemma 5.1]{dend2}, the inclusion $\bext[S]\To \Omega[S]$ is
inner anodyne for any tree $S$ with at least two vertices.
We shall use this to prove that the inclusion
$\Omega[T]_{n-1}\To \Omega[T]_n$ is inner anodyne for $2\leq n\leq N$, which will
prove the proposition.

Let $F_0,\ldots,F_k$ be all subtrees of $T$ having $n$ vertices.
For $0\leq j\leq k$, we put
$$S_j=\bigcup_{0\leq i\leq j}\Omega[F_i]\subset\Omega[T]\, .$$
We shall prove by induction on $j$ that the map
$$S_j\cap\Omega[T]_{n-1}\To S_j$$
is inner anodyne. The case $j=0$ follows from the
identification $\Omega[F_i]\cap\Omega[T]_{n-1}=\bext[F_i]$, $0\leq i\leq k$.
Assume that $j>0$. Note that, since
$$S_{j-1}\cap S_j\cap \Omega[T]_{n-1}=S_{j-1}\cap\Omega[T]_{n-1}\, ,$$
the following diagram is a pushout.
$$\xymatrix{
S_{j-1}\cap\Omega[T]_{n-1}\ar[r]\ar[d]&S_{j-1}\ar[d]\\
S_j\cap\Omega[T]_{n-1}\ar[r]&S_{j-1}\cup(\Omega[F_j]\cap\Omega[T]_{n-1})
}$$
Moreover, since
$\Omega[F_p]\cap\Omega[F_q]\subset\Omega[T]_{n-1}$
for $p\neq q$, we have
$$\Omega[F_j]\cap (S_{j-1}\cup (\Omega[F_j]\cap \Omega[T]_{n-1}))=\Omega[F_j]\cap\Omega[T]_{n-1}\, ,$$
which gives the following pushout square.
$$\xymatrix{
\Omega[F_j]\cap\Omega[T]_{n-1}\ar[r]\ar[d]& \Omega[F_j]\ar[d]\\
S_{j-1}\cup (\Omega[F_j]\cap\Omega[T]_{n-1})\ar[r]& S_j 
}$$
Since the top arrows in these two squares are inner anodyne,
so are the lower ones, and we obtain that the composite
$$S_j\cap\Omega[T]_{n-1}\To S_{j-1}\cup(\Omega[F_j]\cap\Omega[T]_{n-1}) \To S_j$$
is inner anodyne as well. For $j=k$, we conclude that the map
$$\Omega[T]_{n-1}\To S_k\cup \Omega[T]_{n-1}=\Omega[T]_{n}$$
is the pushout of an inner anodyne extension.
\end{proof}
%

\begin{prop}\label{reciproque}
Let $\Wpr$ be a class satisfying the following three conditions.
\begin{itemize}
\item[(i)] The class $\Wpr$ is closed under transfinite
compositions, pushouts and retracts.
\item[(ii)] Any Segal core inclusion belongs to $\Wpr$.
\item[(iii)] For any normal monomorphisms between normal dendroidal
sets $u:X\to Y$ and $v:Y\To Z$, if $u$ and $vu$ are in $\Wpr$, so is $v$.
\end{itemize}
Then any inner anodyne extension belongs to $\Wpr$.
\end{prop}

\begin{proof}
As, by definition, the class of inner anodyne extensions is the smallest
class of maps which satisfies Condition (i) and
contains the inner horn inclusions, it is sufficient
to prove that any inner horn inclusion belongs to $\Wpr$.
For a tree $T$ with at least two vertices, and any internal edge $e$ in $T$, we have
the following natural inclusions.
\begin{equation}\label{reciproque1}
\spine[T]\To \bext[T]\To \Lambda^e[T]\To \Omega[T] 
\end{equation}
We shall prove by induction on the number $|T|$ of vertices of $T$ that
all these inclusions belong to $\Wpr$ (note that this is not so if $T$
has only one vertex). In fact, if $A$ is a union of at least two external
faces of $T$, then $\spine[T]\subset A$, while, if $B$ is the union of
$\bext[T]$ and of a collection of internal faces not including the one
contracting $e$, there are interpolating inclusions
\begin{equation}\label{reciproque2}
\spine[T]\To A \To \bext[T]\To B\To \Lambda^e[T]\To \Omega[T] \, .
\end{equation}
Our induction on $|T|$ will proceed by showing that all these inclusions
belong to $\Wpr$.

To begin with, if $|T|=2$, then $\Omega[T]$ has just two external faces, and
one internal one (given by the edge $e$), so $\spine[T]=A=\bext[T]=B=\Lambda^e[T]$,
whence, as $\Wpr$ contains isomorphisms and Segal core inclusions,
all the inclusions in \eqref{reciproque2} are in $\Wpr$.

Consider now a tree $S$ with $|S|>2$, and assume that, for any tree $T$ such that
$2\leq |T|<|S|$, all the maps in \eqref{reciproque2} are in $\Wpr$. We shall first show
that, for any set $\{R_i\}_{0\leq i\leq j}$ of at least two external faces of $S$, the map
\begin{equation}\label{reciproque3}
\spine[S]\To A=\bigcup_{0\leq i\leq j}\Omega[R_i]
\end{equation}
is in $\Wpr$. For the set of all external faces of $S$, the map \eqref{reciproque3}
is the inclusion $\spine[S]\To \bext[S]$, so that we shall have shown that
this map belongs to $\Wpr$ as well. By Condition (iii),
it then follows that the map $A\To\bext[S]$ also belongs
to $\Wpr$.

To prove that \eqref{reciproque3} is in $\Wpr$, consider the case of
just two distinct external faces $R_0$ and $R_1$ in $S$. Then the map
$$\spine[R_0\cap R_1]\To \Omega[R_0\cap R_1]$$
belongs to $\Wpr$ (also if $R_0\cap R_1$ is a tree with just one vertex),
as does the map $\spine[R_1]\To \Omega[R_1]$.
Now, consider the commutative diagram below.
$$\xymatrix{
\spine[R_0\cap R_1]\ar[r]\ar[d]&\spine[R_1]\ar[dr]\ar[d]&\\
\Omega[R_0]\cap\Omega[R_1]\ar[r]&\spine[R_1]\cup (\Omega[R_0]\cap\Omega[R_1])
\ar[r]&\Omega[R_1]
}$$
The square is a pushout, so the right hand vertical map belongs to $\Wpr$.
As the slanted map also does by assumption, we find that the right hand
horizontal map belongs to $\Wpr$. Next, the two pushout diagrams
\begin{equation*}
\begin{split}
\xymatrix@C=2pt{
\spine[R_1]\cup(\Omega[R_0]\cap\Omega[R_1])\ar[r]\ar[d]&\Omega[R_1]\ar[d]\\
\spine[R_1]\cup \Omega[R_0]\ar[r]& \Omega[R_0]\cup\Omega[R_1]
}\end{split}
\text{and}
\begin{split}
\xymatrix@C=12pt{
\spine[R_0]\ar[r]\ar[d]&\Omega[R_0]\ar[d]\\
\spine[R_0]\cup\spine[R_1]\ar[r]&\Omega[R_0]\cup \spine[R_1]
}\end{split}
\end{equation*}
show that
$$\spine[S]=\spine[R_0]\cup\spine[R_1]\To \Omega[R_0]\cup \Omega[R_1]$$
belongs to $\Wpr$. This shows that \eqref{reciproque3}
belongs to $\Wpr$ for any collection of two distinct external faces.
Consider now a collection $\{R_i\}_{0\leq i\leq j}$ of distinct external faces of $S$,
with $j\geq 1$. We shall prove by induction on $j$ that \eqref{reciproque3}
belongs to $\Wpr$. The case $j=1$ just having been dealt with, we may
assume that $j\geq 2$. By induction, we have that the map
$$\spine[R_0]\To \bigcup_{0<i\leq j}\Omega[R_0\cap R_i]$$
belongs to $\Wpr$, because $R_0\cap R_i$ is an external face of $R_0$
for $0<i\leq j$. Also, the map $\spine[R_0]\To \Omega[R_0]$ is in $\Wpr$, so
that, by Condition (iii), the map
$$ \bigcup_{0<i\leq j}\Omega[R_0\cap R_i]\To \Omega[R_0]$$
belongs to $\Wpr$. By induction, the map
$$\spine[S]\To \bigcup_{0<i\leq j}\Omega[R_i]$$
is in $\Wpr$, whence we deduce from the pushout
$$\xymatrix{
 \bigcup_{0<i\leq j}\Omega[R_0\cap R_i]\ar[r]\ar[d]& \Omega[R_0]\ar[d]\\
 \bigcup_{0<i\leq j}\Omega[R_i]\ar[r]& \bigcup_{0\leq i\leq j}\Omega[R_i]
}$$
that \eqref{reciproque3} is the composite of two maps in $\Wpr$.
This completes the proof of the fact that all the maps of type \eqref{reciproque3}
belong to $\Wpr$.

We now turn to the internal faces of the tree $S$, and let $B$ be the union of
$\bext[S]$ and of a family of internal faces $\Omega[\partial^{a_1}S],\ldots, \Omega[\partial^{a_k}S]$,
given by internal edges $a_1,\ldots,a_k$, all distinct from another internal edge $e$.
We shall prove that $\bext[S]\To B$ belongs to $\Wpr$.
Then the composition $\spine[S]\To\bext[S]\To B$ belongs to $\Wpr$ as well,
and if $B$ contains all internal faces but the one given by $e$,
we find that the inclusion
$\spine[S]\To \Lambda^e[S]$ belongs to $\Wpr$. Since
$\spine[S]\To \Omega[S]$ is in $\Wpr$ by assumption, Condition (iii)
implies that $\Lambda^e[S]\To \Omega[S]$ is in $\Wpr$.
So, to complete the proof of the proposition, it is sufficient to prove that, under
our inductive assumption (that the maps in \eqref{reciproque2}
all belong to $\Wpr$ for smaller trees), the map
\begin{equation}\label{reciproque4}
\bext[S]\To B=\bext[S]\cup \bigcup_{1\leq i\leq k}\Omega[\partial^{a_i}S]
\end{equation}
belongs to $\Wpr$ for any family of internal edges $a_1,\ldots,a_k$ as above.

We proceed by induction on $k$.
If $k=1$, then $S$ has at least two internal edges, hence $\partial^{a_1}S$
has at least two vertices, so, by assumption on trees smaller than $S$,
the map
$$\bext[\partial^{a_1}S]\To \Omega[\partial^{a_1}S]$$
belongs to $\Wpr$.
But
$$\bext[\partial^{a_1}S]=\Omega[\partial^{a_1}S]\cap\bext[S]$$
so we have a pushout diagram
$$\xymatrix{
\bext[\partial^{a_1}S]\ar[r]\ar[d]&\Omega[\partial^{a_1}S]\ar[d]\\
\bext[S]\ar[r]&\bext[S]\cup\Omega[\partial^{a_1}S]\, .
}$$
Thus the inclusion map $\bext[S]\to\bext[S]\cup\Omega[\partial^{a_1}S]$
is in $\Wpr$. This proves the case $k=1$.

If $k>1$, then we have
$$\left(\bext[S]\cup \bigcup_{1\leq i<k}\Omega[\partial^{a_i}S]\right)\cap\Omega[\partial^{a_{k}}S]=
\bext[\partial^{a_k}S]\cup \bigcup_{1\leq i\leq k}\Omega[\partial^{a_i}\partial^{a_k}S]\, ,$$
so the diagram
$$\xymatrix{
\bext[\partial^{a_k}S]\cup \bigcup_{1\leq i\leq k}\Omega[\partial^{a_i}\partial^{a_k}S]
\ar[r]\ar[d]&\Omega[\partial^{a_k}S]\ar[d]\\
\bext[S]\cup \bigcup_{1\leq i<k}\Omega[\partial^{a_i}S]\ar[r]&
\bext[S]\cup \bigcup_{1\leq i\leq k}\Omega[\partial^{a_i}S]
}$$
is a pushout. Moreover, the family $\{\partial^{a_i}\partial^{a_k}S\}_{1\leq i<k}$
of internal faces of $T=\partial^{a_k}S$ misses the edge $e$, so that, by assumption, all the
inclusions in \eqref{reciproque2} belong to $\Wpr$.
We conclude that the top arrow in the pushout above belongs to $\Wpr$, whence
so does the bottom arrow. By induction on $k$, the map
$$\bext[S]\To\bext[S]\cup \bigcup_{1\leq i< k}\Omega[\partial^{a_i}S]$$
is in $\Wpr$, and as $\Wpr$ is closed under composition,
we find that \eqref{reciproque4} is in $\Wpr$.
\end{proof}

\begin{cor}\label{Segalchardendnerve}
A dendroidal set $X$ is the nerve of an operad if and only if, for any tree $T$, the map
$$X_T=\Hom_\dset(\Omega[T],X)\To \Hom_\dset(\spine[T],X)$$
is bijective.
\end{cor}

\begin{proof}
Assume that $X=\mathit{N}_d(P)$ for a (coloured) operad $P$.
Then, as the functor $\tau_d$ sends inner anodyne inclusions to isomorphisms,
it follows from Proposition \ref{segalcoreanodyne} that the map
$\Hom_\dset(\Omega[T],X)\To \Hom_\dset(\spine[T],X)$ is bijective.
The converse follows easily from Proposition \ref{reciproque}
and from the characterization of dendroidal nerves
given by \cite[Proposition 5.3 and Theorem 6.1]{dend2}.
\end{proof}

\begin{rem}
The preceding corollary
gives in particular the well known characterization of small categories
as simplicial sets satisfying the Grothendieck-Segal condition:
given a simplicial set $X$, then, for $T=n$
with $n\geq 1$, we have
$$\Hom_\dset(\spine[T],X)=\underset{\text{$n$ times}}{
\underbrace{X_1\times_{X_0}X_1\times_{X_0}\dots \times_{X_0}X_1}}\, .$$
\end{rem}

\section{Equivalences of $\infty$-operads}

\begin{paragr}\label{reminderHoms}
For an $\infty$-category $X$, we denote by $k(X)$ the maximal
Kan complex contained in $X$; see \cite[Corollary 1.5]{joyal}.
Recall that, if $A$ and $X$ are two dendroidal sets, $\sHom(A,X)$ denotes their
internal Hom (with respect to the Boardman-Vogt tensor product of
dendroidal sets).

Given two dendroidal sets $A$ and $X$, we write
$$\map(A,X)=i^*\sHom(A,X)\, .$$
If $X$ is an $\infty$-operad, and if $A$
is normal, then $\sHom(A,X)$ is an $\infty$-operad, so that $\map(A,X)$ is
an $\infty$-category; see \cite[Theorem 9.1]{dend2}.

For an $\infty$-operad $X$ and a simplicial set $K$, we will write
$X^{(K)}$ for the subcomplex of $\sHom(i_!(K),X)$ which consists
of dendrices
$$a:\Omega[T]\times i_!(K)\To X$$
such that, for any edge $u$ in the tree $T$, the induced map
$$a_u:K\To i^*(X)$$
factors through $k(i^*(X))$ (i.e. all the $1$-cells in the image
of $a_u$ are weakly invertible in $i^*(X)$).

For an $\infty$-operad $X$ and a normal dendroidal set $A$, we
write $k(A,X)$ for the subcomplex of $\map(A,X)$ which consists
of maps
$$u: A\otimes i_!(\Delta[n])\To X$$
such that, for all vertices $a$ of $A$ (i.e. maps $a:\eta\To A$), the induced map
$$u_a:\Delta[n]\To i^*(X)$$
factors through $k(i^*(X))$.
So, by definition, for any normal dendroidal set $A$, any simplicial set $K$,
and any $\infty$-operad $X$, there is a natural bijection:
\begin{equation}\label{dendanalogofk1}
\Hom_{\sset}(K,k(A,X))\simeq\Hom_{\dset}(A,X^{(K)})\, .
\end{equation}
Furthermore, by virtue of \cite[6.8]{dend3}, we have the equality:
\begin{equation}\label{dendanalogofk2}
k(A,X)=k(\map(A,X))
\end{equation}
(in particular, $k(A,X)$ is a Kan complex).
\end{paragr}

\begin{paragr}\label{paragr:recallkAXframe}
Recall that, in any model category $\C$,  given a cofibrant object $A$ and a fibrant
object $X$, one of the models of the mapping space $\mathit{Map}(A,X)$
is the simplicial set defined by
\begin{equation}\label{eq:Map1}
\mathit{Map}(A,X)_n=\Hom_\C(A,X_n)\, ,
\end{equation}
where $X_\bullet$ is a Reedy fibrant resolution of $X$ (where $X$ is seen as
a simplicially constant object of $\C^{\op{\Delta}}$); see \cite[5.4.7 and 5.4.9]{Ho},
for instance. In the case where $\C=\dset$ and $X$ is an $\infty$-operad, then
the simplicial dendroidal set $X^{(\smp{\bullet})}$ is Reedy fibrant, and, for any
integer $n\geq 0$, the map $X\To X^{(\smp{n})}$ is a weak equivalence
(this follows immediately from \cite[Corollary 6.9]{dend3}).
In particular, we get:
\end{paragr}

\begin{prop}\label{kAXMapAX}
If $A$ is a normal dendroidal set and $X$ an $\infty$-operad, then
there is a natural isomorphism (in fact, identity) of simplicial sets
$$k(A,X)\simeq \mathit{Map}(A,X)\, .$$
\end{prop}

\begin{proof}
This follows immediately from the identifications \eqref{dendanalogofk1}
and \eqref{eq:Map1}.
\end{proof}

\begin{lemma}\label{kderinthom}
Let $X$ be an $\infty$-operad. 
\begin{itemize}
\item[(i)] For any simplicial set $K$, and for any normal dendroidal set
$A$, there is a natural bijection
$$\Hom_{\ho(\sset)}(K,k(A,X))\simeq \Hom_{\ho(\dset)}(A,X^{(K)})\, .$$
\item[(ii)] For any cofibration (resp. trivial cofibration)
between normal dendroidal sets $A\To B$, the map
$k(B,X)\To k(A,X)$ is a fibration (resp. a trivial fibration)
between Kan complexes.
\item[(iii)] For any pushout of normal dendroidal sets
$$\xymatrix{
A\ar[r]\ar[d]_i& A'\ar[d]^{i'}\\
B\ar[r]& B'
}$$
with $i$ a cofibration, the commutative square
$$\xymatrix{
k(B',X)\ar[r]\ar[d]& k(A',X)\ar[d]\\
k(B,X)\ar[r]& k(A,X)
}$$
is a pullback.
\item[(iv)] For any sequence of cofibrations between normal
dendroidal sets
$$A_0\To A_1\To \cdots\To A_n\To A_{n+1}\To \cdots \, ,$$
the map
$$k(\varinjlim_n A_n,X)\To \varprojlim_n k(A_n,X)$$
is an isomorphism.
\end{itemize}
\end{lemma}

\begin{proof}
This follows immediately from Proposition \ref{kAXMapAX}, using the
general properties of mapping space functors; see for instance
\cite[Proposition 5.4.1 and Theorem 5.4.9]{Ho}.
\end{proof}

\begin{thm}\label{charactwoperequiv1}
Let $f:X\To Y$ be a morphism between $\infty$-operads.
The following conditions are equivalent.
\begin{itemize}
\item[(a)] For any integer $n\geq 0$, the map $k(\Omega[C_n],X)\To k(\Omega[C_n],Y)$
is a simplicial homotopy equivalence, as is the map $k(\eta,X)\To k(\eta,Y)$.
\item[(b)] For any tree $T$, the map $k(\Omega[T],X)\To k(\Omega[T],Y)$
is a simplicial homotopy equivalence.
\item[(c)] For any normal dendroidal set $A$, the map $k(A,X)\To k(A,Y)$
is a simplicial homotopy equivalence.
\item[(d)] The map $f:X\To Y$ is a weak operadic equivalence.
\end{itemize}
\end{thm}

\begin{proof}
Assume condition (a). We claim that, for any tree $T$, the induced map
$k(\spine[T],X)\To k(\spine[T],Y)$ is a simplicial homotopy
equivalence. Note that, for $T=[0]$, this is a special case of (a). Therefore,
to prove this, we may assume that $T$ has at least one vertex.
Let $v_1,\dots,v_k$ the vertices of $T$, and, for $1\leq i\leq k$,
write $n_i$ for the number of input edges of $v_i$ in $T$. We then have
$$\spine[T]=\bigcup_{1\leq i\leq k}\Omega[C_{n_i}]\, .$$
Moreover, for two indices $i\neq j$, the intersection
$\Omega[C_{n_i}]\cap \Omega[C_{n_j}]$ is either empty
or isomorphic to $\eta$. For $1\leq j\leq k$, define
$$A_i=\bigcup_{1\leq i\leq j}\Omega[C_{n_i}]\, .$$
For $1< i\leq k$, there is a pushout square
$$\xymatrix{
A_{i-1}\cap  \Omega[C_{n_i}]\ar[r]\ar[d] &  \Omega[C_{n_i}]\ar[d]\\
A_{i-1}\ar[r] & A_i
}$$
in which the intersection $ A_{i-1} \cap \Omega[C_{n_i}]$ is
isomorphic to a finite sum of $\eta$'s.
By the cube lemma (see the dual version of \cite[Lemma 5.2.6]{Ho}),
using properties (ii) and (iii) of Lemma \ref{kderinthom}, we obtain by induction
on $i$ that the maps
$$k(A_i,X)\To k(A_i,Y)$$
are simplicial homotopy equivalences.
In particular, for $i=k$, this means that the map
$$k(\spine[T],X)\To k(\spine[T],Y)$$
is a simplicial homotopy equivalence.
Thus, since the vertical maps in the commutative square
$$\xymatrix{
k(\Omega[T],X)\ar[r]\ar[d]& k(\Omega[T],Y)\ar[d]\\
k(\spine[T],X)\ar[r]&k(\spine[T],Y) 
}$$
are simplicial homotopy equivalences as well
(by Proposition \ref{segalcoreanodyne} and Lemma \ref{kderinthom}~(ii)),
this proves (b).

The fact that condition (b) implies condition (c) follows by
similar arguments\footnote{Another proof consists to
see that, by \cite[1.7]{dend3},
for a normal dendroidal set $A$, the category $\Omega/A$
is a regular skeletal category in the sense of \cite[8.2.3]{Ci3}, from which
we deduce that $A$ is the homotopy colimit of the
$\Omega[T]$'s over $A$ (see \cite[8.2.9]{Ci3}), and we can use
Proposition \ref{kAXMapAX} to see that the functor $k(-,X)$
turns homotopy colimits into homotopy limits (see for instance
\cite[6.13]{Ci1}), which implies the result.}
from Lemma \ref{kderinthom}, using the skeletal filtration
of normal dendroidal sets~\cite[section 4]{dend2}.

A reformulation of Lemma \ref{kderinthom}~(i) in the particular case
where $K=\Delta[0]$ is that,
for a normal dendroidal set $A$ and an $\infty$-operad $X$,
we have a natural bijection
$$\pi_0(k(A,X))\simeq\Hom_{\ho(\dset)}(A,X)\, .$$
It thus follows from the Yoneda lemma that condition (c) implies condition (d).

Finally, the fact that condition (d) implies condition (c) (and, therefore, condition (a)) is obvious:
it follows from Lemma \ref{kderinthom}~(ii) and from Ken~Brown's Lemma \cite[Lemma 1.1.12]{Ho}
that the functor $k(A,-)$ sends weak operadic equivalences between $\infty$-operads
to simplicial homotopy equivalences.
\end{proof}

\begin{paragr}\label{def:Homininftyoperads}
Let $X$ be an $\infty$-operad.
Given an $(n+1)$-tuple of $0$-cells $(x_1,\ldots,x_n,x)$ in $X$,
the space of maps $X(x_1,\ldots,x_n;x)$ is obtained
by the pullback below, in which the map $p$ is
the map induced by the inclusion $\eta\amalg\cdots\amalg\eta\To\Omega[C_n]$
(with $n+1$ copies of $\eta$, corresponding to the $n+1$ objects
$(a_1,\ldots,a_n,a)$ of $C_n$; see \eqref{corolla}).
$$\xymatrix{
X(x_1,\ldots,x_n;x)\ar[r]\ar[d]&\sHom(\Omega[C_n],X)\ar[d]^p\\
\eta\ar[r]_(.45){(x_1,\ldots,x_n,x)}&X^{n+1}
}$$
Using the identification $\sset=\dset/\eta$, we shall consider
$X(x_1,\ldots,x_n;x)$ as a simplicial set.
Observe that $X(x_1,\ldots,x_n;x)$ is actually a Kan complex
(see \cite[Proposition 6.13]{dend3}).
\end{paragr}

\begin{definition}\label{def:equivinftyoperd}
Let $f:X\To Y$ be a morphism of $\infty$-operads.

The map $f$ is \emph{fully faithful} if, for any
$(n+1)$-tuple of $0$-cells $(x_1,\ldots,x_n,x)$ in $X$,
the morphism
$$X(x_1,\ldots,x_n;x)\To Y(f(x_1),\ldots,f(x_n);f(x))$$
is a simplicial homotopy equivalence.

The map $f$ is \emph{essentially surjective} if
the functor underlying the morphism of operads $\tau_d(X)\To \tau_d(Y)$ is essentially surjective.
\end{definition}

\begin{rem}\label{rempinotofk}
For an $\infty$-operad $X$, the set of isomorphism classes of objects in
$\tau_d(X)$ is in bijection with the set $\pi_0(k(i^*X))$; see \cite[4.1]{dend3}.
The condition of essential surjectivity is thus equivalent to the fact that
the map $k(i^*X))\To k(i^*Y))$ induces a surjection on connected components.

Moreover, by virtue of \cite[Proposition 6.14]{dend3}, we have natural bijections
$$\pi_0(X(x_1,\ldots,x_n;x))\simeq \tau_d(X)(x_1,\ldots,x_n;x)\, .$$
As a consequence, if $f:X\To Y$ is fully faithul, so is the induced morphism of operads
$\tau_d(f):\tau_d(X)\to \tau_d((Y)$. Therefore, if $f:X\To Y$ is fully faithful
and essentially surjective, then the map $k(i^*X))\To k(i^*Y))$
induces a bijection on connected components.
\end{rem}

We recall the following well known fact:

\begin{lemma}\label{criteriumweakequi0}
A commutative square of simplicial sets
$$\xymatrix{
X\ar[r]^u\ar[d]_p&X'\ar[d]^{p'}\\
Y\ar[r]_v&Y'
}$$
in which $p$ and $p'$ are Kan fibrations is a homotopy
pullback square if and only if, for
any $0$-simplex $y$ of $Y$, the map between the corresponding
fibers
$$p^{-1}(y)=X_y\To X'_{v(y)}=p^{\prime-1}(v(y))$$
is a simplicial homotopy equivalence.
\end{lemma}

A direct consequence of the preceding lemma is:

\begin{lemma}\label{criteriumweakequi}
Consider a commutative square of simplicial sets.
$$\xymatrix{
X\ar[r]^u\ar[d]_p&X'\ar[d]^{p'}\\
Y\ar[r]_v&Y'
}$$
Assume furthermore that $p$ and $p'$ are Kan fibrations, and that
the map $v$ is a weak homotopy equivalence.
Then the map $u$ is a weak homotopy equivalence if and only if, for
any $0$-simplex $y$ of $Y$, the map between the corresponding
fibers
$$X_y\To X'_{v(y)}$$
is a simplicial homotopy equivalence.
\end{lemma}

\begin{thm}\label{charactwoperequiv2}
Let $f:X\To Y$ be a morphism between two $\infty$-operads.
Then $f$ is a weak operadic equivalence if and only if
it is fully faithful and essentially surjective.
\end{thm}

\begin{proof}
Given an $\infty$-operad $X$, and an $(n+1)$-tuple $(x_1,\ldots,x_n,x)$ of objects of $X$,
i.e. a $0$-simplex of $k(i^*X)^{n+1}=k(\bord\Omega[C_n],X)$,
we have the following diagram in which the right hand square
is a pullback square (see \cite[Remark 6.2 and Corollary 6.8]{dend3}).
$$\xymatrix{
X(x_1,\ldots,x_n;x)\ar[r]\ar[d]&k(\Omega[C_n],X)\ar[r]\ar[d]&\map(\Omega[C_n],X)\ar[d]\\
\eta\ar[r]_{(x_1,\ldots,x_n,x)}&k(i^*X)^{n+1}\ar[r]&i^*X^{n+1}
}$$
Hence the left hand square above is a pullback as well.
Assume that $f$ is fully faithful and essentially surjective. We will first prove that the induced morphism
$k(i^*X)\To k(i^*Y)$ is a simplicial homotopy equivalence.
As the corresponding map $\pi_0(k(i^*X))\To\pi_0(k(i^*Y))$ is bijective (see \ref{rempinotofk}),
it is sufficient to prove that, for any $0$-simplex $x$ of $k(i^*X)$, the map of loop spaces
$$\Omega(k(i^*X),x)\To\Omega(k(i^*Y),f(x))$$
is a weak homotopy equivalence.
For this purpose, it will be sufficient to prove that the commutative
square
\begin{equation}\label{charactwoperequiv2proof1}\begin{split}
\xymatrix{
\Omega(k(i^*X),x)\ar[r]\ar[d]&\Omega(k(i^*Y),f(x))\ar[d]\\
X(x;x)\ar[r]&Y(f(x);f(x))
}\end{split}\end{equation}
is homotopy cartesian. In general,
the set of connected components of the Kan complex $X(x_1,\ldots,x_n;x)$
is in bijection with the set $\tau_d(X)(x_1,\ldots,x_n;x)$; see
\cite[Proposition 6.14]{dend3}.
We can thus describe the loop space $\Omega(k(i^*X),x)$ as the disjoint union
of the connected components of $X(x;x)$ which correspond to
automorphisms of $x$ in the category underlying the operad $\tau_d(X)$.
In particular, the inclusion $\Omega(k(i^*X),x)\subset X(x;x)$ is a Kan
fibration between Kan complexes. Using the fact that the functor underlying the map
$\tau_d(X)\To\tau_d(Y)$ is full faithful (whence conservative),
we see that the map
$$\pi_1(k(i^*X),x)=\pi_0(\Omega(k(i^*X),x))\To \pi_0(\Omega(k(i^*Y),f(x)))=\pi_1(k(i^*Y),f(x))$$
is bijective. The square \eqref{charactwoperequiv2proof1} is thus cartesian, and,
as its vertical maps are Kan fibrations, it is homotopy cartesian as well.
Therefore, the map $k(i^*X)\To k(i^*Y)$ is a simplicial homotopy equivalence.
As a consequence, for any corolla $C_n$, we also have simplicial
homotopy equivalences
$$k(i^*X)^{n+1}=k(\bord\Omega[C_n],X)\To k(i^*Y)^{n+1}=k(\bord\Omega[C_n],Y)\, .$$
By applying Lemma \ref{criteriumweakequi} to the commutative squares
\begin{equation}\label{charactwoperequiv2proof2}\begin{split}\xymatrix{
k(\Omega[C_n],X)\ar[r]\ar[d]&k(\Omega[C_n],Y)\ar[d]\\
k(\bord\Omega[C_n],X)\ar[r]&k(\bord\Omega[C_n],Y)
}\end{split}\end{equation}
we conclude that the maps $k(\Omega[C_n],X)\To k(\Omega[C_n],Y)$
are all simplicial homotopy equivalences. The characterization
given by condition (a) of Theorem \ref{charactwoperequiv1} thus implies that
$f$ is a weak operadic equivalence.

For the converse, we just apply again Lemma \ref{criteriumweakequi}
to the commutative squares \eqref{charactwoperequiv2proof2}.
\end{proof}

\begin{rem}
As we saw implicitly in the proof above, Theorem \ref{charactwoperequiv1}
and Lemma \ref{criteriumweakequi0} lead to a characterization
of fully faithful maps: a morphism between $\infty$-operads $X\To Y$
is fully faithful if and only if the commutative squares of shape \eqref{charactwoperequiv2proof2}
are homotopy pullback squares of Kan complexes for any $n\geq 0$.
\end{rem}

\section{Locally constant simplicial dendroidal sets}\label{section:locconstant}

\begin{paragr}
Let $\sdset=\dset^{\op{\cats}}\simeq\pref{\cats\times\Omega}$
be the category of simplicial dendroidal sets.
The category $\dset$ (resp. $\sset$) of dendroidal (resp. simplicial) sets is naturally
embedded in $\sdset$, by viewing a dendroidal (resp. simplicial) set as a constant
simplicial (resp. dendroidal) object in $\dset$ (resp. in $\sset$).
For a simplicial dendroidal set $X$, a tree $T$, and an integer $n\geq 0$, the evaluation
of $X$ at $(T,n)$ will often be denoted by $X(T)_n$.

Given a simplicial dendroidal set $X$, we denote by
$$\op{\sset}\To\dset\ , \quad K \longmapsto X^K$$
the (essentially) unique limit preserving functor
which sends $\smp{n}$ to $X^{\smp{n}}:=X_n$.

Starting from the model category structure on $\dset$, and using that
the category of simplices $\cats$ is a Reedy category,
one obtains the Reedy model structure on $\sdset$; see \cite[Theorem 5.2.5]{Ho}.
We shall call this structure the \emph{simplicial Reedy model
category structure on $\sdset$.}
By definition, the weak equivalences are the termwise weak operadic equivalences
(by evaluating at simplices), while the fibrations (resp. the trivial
fibrations) are the maps $X\To Y$ such that, for any integer $n\geq 0$, the map
$$X^{\smp{n}}\To X^{\partial\smp{n}}\times_{Y^{\partial\smp{n}}}Y^{\smp{n}}$$
is a fibration (resp. a trivial fibration) in $\dset$. In other words, we have:
\end{paragr}

\begin{prop}\label{generators}
The simplicial Reedy model structure on $\sdset$ is a cofibrantly generated
model category. A generating set of cofibrations of $\sdset$ is given by the inclusions
$$\partial\smp{n}\times\Omega[T]\cup \smp{n}\times\partial\Omega[T]\To
\smp{n}\times\Omega[T]$$
for any integer $n\geq 0$ and any tree $T$, while, if $K$ is a generating set
of trivial cofibrations (between normal dendroidal sets) in $\dset$, then a
generating set of trivial cofibrations of $\sdset$ is given by the inclusions
$$\partial\smp{n}\times B\cup \smp{n}\times A\To
\smp{n}\times B$$
for any integer $n\geq 0$ and any map $A\To B$ in $K$.
\end{prop}

\begin{cor}\label{simplnormalmono}
The cofibrations of the simplicial Reedy model category structure on $\sdset$
are the termwise normal monomorphisms of dendroidal sets.
\end{cor}

\begin{proof}
It is well known that any Reedy cofibration is a termwise
cofibration. Therefore, it is sufficient to prove
that any termwise normal monomorphism is a cofibration
of the simplicial Reedy model category structure on $\sdset$.
As both $\Omega$ and $\Delta$ are skeletal categories in the sense
of \cite[8.1.1]{Ci3}, so is there product \cite[8.1.7]{Ci3}.
Moreover, for any integer $n\geq 0$ and any tree $T$, the boundary of
the representable presheaf $\smp{n}\times\Omega[T]$ is nothing but the
presheaf $\partial\smp{n}\times\Omega[T]\cup \smp{n}\times\partial\Omega[T]$.
Therefore, the Reedy cofibrations of $\sdset$ are the normal
monomorphisms in the absolute sense (see \cite[8.1.30 and 8.1.35]{Ci3} for $A=\cats\times\Omega$).
Thus Reedy cofibrations are precisely the monomorphisms $X\to Y$ in
$\sdset$ such that, for any integer $n\geq 0$ and any tree $T$, any
non-degenerate element $y\in Y_{n,T}$ which does not belong to the image of $X_{n,T}$ has a trivial
stabilizer in $\mathrm{Aut}(([n],T))=\mathrm{Aut}(T)$.
On the other hand, a monomorphism of simplicial dendroidal sets $X\to Y$
is termwise normal if and only if, for any integer $n\geq 0$
and any tree $T$, any element $y\in Y_{n,T}$ which does not belong to the image
of $X_{n,T}$ has a trivial stabilizer.
Therefore, any termwise normal monomorphism is a Reedy cofibration.
\end{proof}

\begin{paragr}
In the sequel, we shall simply call \emph{normal monomorphisms} the
cofibrations of the simplicial Reedy model category structure on $\sdset$.
\end{paragr}

\begin{rem}\label{Reedyfibrantoperad}
Any fibrant object of the simplicial Reedy model category on $\sdset$ is
termwise fibrant (\emph{i.e.} is termwise an $\infty$-operad).
\end{rem}

\begin{definition}\label{def:loccstmc}
We define the \emph{locally constant model category structure} on $\sdset$
as the left Bousfield localization of the simplicial Reedy model category structure on $\sdset$
by the set of projections $\smp{n}\times\Omega[T]\To \Omega[T]$, for any tree $T$ and
any integer $n\geq 0$ (see \cite{Hir} for the general theory of left Bousfield localization
of model categories).
\end{definition}

\begin{prop}\label{locconstReedyfibrant}
Let $X$ a simplicial dendroidal set. Assume that $X$ is
fibrant for the simplicial Reedy model category structure.
Then the following conditions are equivalent:
\begin{itemize}
\item[(i)] the map from $X$ to the terminal object has the right lifting property
with respect to the inclusions
$$\Lambda^k[n] \times\Omega[T]\cup \smp{n}\times\partial\Omega[T]\To
\smp{n}\times\Omega[T]$$
for any tree $T$ and for any integers $n\geq 1$ and $0\leq k\leq n$;
\item[(ii)] for any integer $n\geq 0$, the map $X_0\To X_n$ is an equivalence
of $\infty$-operads;
\item[(iii)] $X$ is fibrant for the locally constant model category structure on $\sdset$.
\end{itemize}
\end{prop}

\begin{proof}
The equivalence between conditions (ii) and (iii) holds by definition
of left Bousfield localizations, as, for any tree $T$,
we have a natural identification in $\ho(\sset)$
$$\mathit{Map}(\Omega[T],X_n)=\mathit{Map}(\smp{n}\times \Omega[T],X)$$
for any integer $n\geq 0$.

Next, it follows easily from \cite[Corollary 2.1.20]{Ci3} and \cite[Lemma 4.2.4]{Ho} that
the class of monomorphisms $K\To L$ in $\sset$ such that, for any tree $T$, the map
$$K \times\Omega[T]\cup L \times\partial\Omega[T]\To
L \times\Omega[T]$$
is a trivial cofibration of the locally constant model category structure on $\sdset$,
contains the class of trivial cofibrations of the usual model category structure
on $\sset$. Therefore, condition (iii) implies condition (i).
Conversely, as the horn inclusions generate the trivial cofibrations of the usual model
category structure on $\sset$, it is clear that condition (i) implies
condition (iii).
\end{proof}

\begin{prop}\label{trivQuillenequivloccst}
The inclusion $\dset\subset \sdset$ is a left Quillen equivalence
from the model category for $\infty$-operads to the
locally constant model category structure.
Moreover, this inclusion preserves and detects weak equivalences
between arbitrary objects.
\end{prop}

\begin{proof}
This inclusion functor is left adjoint to the
evaluation at zero functor
$$\mathit{ev}_0:\sdset\to\dset\ , \quad X\longmapsto X_0\, .$$
Note first that the inclusion functor $\dset\subset \sdset$ is a left Quillen
functor which preserves weak equivalences:
by virtue of Corollary \ref{simplnormalmono}, it preserves cofibrations,
and it preserves weak equivalences by definition of the
locally constant model structure.
Thus, to finish the proof, it is sufficient to check the following two
properties:
\begin{itemize}
\item[(a)] for any fibrant object $X$ of
the locally constant model structure, the natural map
$X_0\to X$ is a weak equivalence;
\item[(b)] for any fibrant object $X$ in $\dset$, there exists a weak
equivalence $X\to Y$ in $\sdset$ with $Y$ fibrant in
the locally constant model structure, such that
the induced map $X\to Y_0$ is a weak equivalence of
dendroidal sets.
\end{itemize}
Property (a) follows from the characterization of fibrant objects
given by condition (ii) of the previous proposition.
Property (b) is a particular instance of the existence of Reedy fibrant
resolutions\footnote{One can also
prove (b) more explicitely from \cite[Corollary 6.9]{dend3}:
given an $\infty$-operad $X$,
we can consider the simplicial dendroidal set $Y$ defined by
$Y_{n,T}=k(\Omega[T],X)_n$ for any integer $n\geq 0$ and any
tree $T$ (see Paragraph \ref{paragr:recallkAXframe} above).}.
\end{proof}
%

\section{Dendroidal Segal spaces}\label{section:dendsegalspaces}

\begin{paragr}\label{paragr:defexpA}
We shall now consider different model category structures on the category
of simplicial dendroidal sets.

Given a simplicial dendroidal set $X$, let us denote by
$$\op{\dset}\To \sset \ , \quad A\longmapsto X^A$$
the (essentially) unique limit preserving functor which sends a tree $T$ to
$X^{\Omega[T]}:=X(T)$.

Starting from the usual model category structure on the
category of simplicial sets, we first have:
\end{paragr}

\begin{prop}\label{dendroidalReedy}
The category $\sdset$ admits a cofibrantly generated
and proper model category structure whose weak equivalences
are the termwise simplicial weak homotopy equivalences (\emph{i.e.} the maps
$X\To Y$ such that, for any tree $T$, the map $X_T\To Y_T$ is a simplicial weak
homotopy equivalence), and whose cofibrations are the normal monomorphisms.
Moreover, a morphism of simplicial dendroidal sets $X\To Y$ is a fibration
(resp. a trivial fibration) if and only if, for any tree $T$, the map
\begin{equation}\label{eq:genReedyfib}
X^{\Omega[T]}\To X^{\partial\Omega[T]}\times_{Y^{\partial\Omega[T]}} Y^{\Omega[T]}
\end{equation}
is a Kan fibration (resp. a trivial Kan fibration). In other words, a set of generators
for cofibrations is provided by the maps
$$\partial\smp{n} \times\Omega[T]\cup \smp{n}\times\partial\Omega[T]\To
\smp{n}\times\Omega[T]$$
for any tree $T$ and for any integer $n\geq 0$, while a generating set of
trivial cofibrations is given by the maps
$$\Lambda^k[n] \times\Omega[T]\cup \smp{n}\times\partial\Omega[T]\To
\smp{n}\times\Omega[T]$$
for any tree $T$ and for any integers $n\geq 1$ and $0\leq k\leq n$.
\end{prop}

\begin{proof}
As the category $\Omega$ is a \emph{generalized Reedy category}
(see \cite[Example 2.8]{BMReedy}), the model category above
can be obtained as special case of \cite[Theorem 1.6]{BMReedy}.
%
\end{proof}

\begin{rem}
We shall call the model category of Proposition \ref{dendroidalReedy}
the \emph{generalized Reedy model category structure} on $\sdset$.
\end{rem}

\begin{definition}\label{defSegal}
We define the \emph{model category structure for dendroidal Segal spaces}
as the left Bousfield localization of the generalized Reedy model category
on $\sdset$ by the set of maps $\spine[T]\To \Omega[T]$, for any tree $T$
(we consider $\dset$ as full subcategory of $\sdset$ via the
obvious inclusion $\dset\subset\sdset$).
A fibrant object for this model category will be called a \emph{dendroidal Segal space}.
\end{definition}

\begin{prop}\label{Segalspinesinnerhorns}
The model category structure for dendroidal Segal spaces
is the left Bousfield localization of the generalized Reedy model category
structure on $\sdset$ by the set of maps $\Lambda^e[T]\To \Omega[T]$,
for any tree $T$ with an inner edge $e$.
\end{prop}

\begin{proof}
It follows immediately from Proposition \ref{reciproque}
that, for any tree $T$ with given inner edge $e$,
the map $\Lambda^e[T]\To \Omega[T]$
is a weak equivalence of the model category structure for dendroidal Segal spaces.
Conversely, let $\Wpr$ be the class of maps of dendroidal sets
which are weak equivalences in the left Bousfield localization of the generalized
Reedy model category on $\sdset$ by the set of maps $\Lambda^e[T]\To \Omega[T]$,
for any tree $T$ with an inner edge $e$.
It is clear that any inner anodyne extension in $\dset$ belongs to $\Wpr$,
so that, by virtue of Proposition \ref{segalcoreanodyne}, for any tree $T$, the
inclusion $\spine[T]\To\Omega[T]$ is in $\Wpr$.
\end{proof}

\begin{cor}\label{chardendSegal}
Let $X$ be a simplicial dendroidal set. Assume that $X$ is
fibrant for the generalized Reedy model category structure.
Then the following conditions are equivalent:
\begin{itemize}
\item[(i)] $X$ is a dendroidal Segal space;
\item[(ii)] for any tree $T$, the map $X^{\Omega[T]}\To X^{\spine[T]}$
is a trivial Kan fibration in $\sset$;
\item[(iii)] for any tree $T$ with a given inner edge $e$, the map
$X^{\Omega[T]}\To X^{\Lambda^e[T]}$ is a trivial Kan fibration in $\sset$.
\end{itemize}
\end{cor}

\begin{proof}
For any normal dendroidal set $A$, there is a canonical identification
$$\mathit{Map}(A,X)=X^A\, .$$
The corollary thus follows from the definition
of left Bousfield localizations and from Proposition \ref{Segalspinesinnerhorns}.
\end{proof}

\begin{prop}\label{charweakequivSegalspaces}
A morphism of dendroidal Segal spaces $X\To Y$ is a weak
equivalence if and only if, its evaluation at $T$
$$X(T)=X^{\Omega[T]}\To Y^{\Omega[T]}=Y(Y)$$
is a simplicial homotopy equivalence for $T=\eta$
as well as for $T=C_n$, $n\geq 0$.
\end{prop}

\begin{proof}
A morphism of dendroidal Segal spaces is a weak equivalence of
the model structure for dendroidal Segal spaces if and only if
it is a weak equivalence of the generalized Reedy model category
structure. In other words, a morphism of dendroidal Segal spaces $X\To Y$ is a weak
equivalence if and only if, its evaluation at $T$
$$X^{\Omega[T]}\To Y^{\Omega[T]}$$
is a simplicial homotopy equivalence for any tree $T$.
By virtue of condition (ii) of the preceding corollary, we see that
evaluating a dendroidal Segal space at a tree $T$ gives the same information
as evaluating at $\spine[T]$. We easily conclude the proof from the fact that
$\spine[T]$ is a (homotopy) colimit of dendroidal sets of shape $\eta$ or $\Omega[C_n]$, $n\geq 0$.
\end{proof}

\begin{paragr}
If $X$ is a dendroidal Segal space, and if $(x_1,\ldots,x_n,x)$ is an $(n+1)$-tuple of
elements of $X(\eta)_0^{n+1}$, we define $X(x_1,\ldots,x_n;x)$ by the following
pullback
$$\xymatrix{
X(x_1,\ldots,x_n;x)\ar[r]\ar[d]&X(C_n)\ar[d]\\
\Delta[0]\ar[r]_{(x_1,\ldots,x_n,x)}&X(\eta)^{n+1}
}$$
in which the map $X(C_n)\To X(\eta)^{n+1}=X^{\bord \Omega[C_n]}$ is the map
induced by the inclusion $\coprod_{n+1}\eta= \bord\Omega[C_n]\To \Omega[C_n]$.
As $X(C_n)\To X(\eta)^{n+1}$ is a Kan fibration, the pullback square above is
homotopy cartesian, and $X(x_1,\ldots,x_n;x)$ is a Kan complex.
\end{paragr}

\begin{definition}\label{def:fullyfaithfulSegalspaces}
A morphism of dendroidal Segal spaces $f:X\To Y$ is \emph{fully faithful} if, for any
$(n+1)$-tuple of $0$-cells $(x_1,\ldots,x_n,x)$ in $X(\eta)$,
the morphism
$$X(x_1,\ldots,x_n;x)\To Y(f(x_1),\ldots,f(x_n);f(x))$$
is a simplicial homotopy equivalence.

A morphism of dendroidal Segal spaces is \emph{a weak equivalence on objects}
if its evaluation at $\eta$ is a simplicial weak equivalence.
\end{definition}

\begin{cor}\label{cor:charweakequivSegalspaces}
A morphism of dendroidal Segal spaces is a weak equivalence if and only
if it is fully faithful as well as a weak equivalence on objects.
\end{cor}

\begin{proof}
This follows immediately from Proposition \ref{charweakequivSegalspaces}
and from Lemma \ref{criteriumweakequi}.
\end{proof}

\section{Complete dendroidal Segal spaces}\label{dendrezk}

\begin{paragr}
Recall the dendroidal interval $J_d=i_!(J)$,
where $J= N\pi_1(\Delta[1])$ denotes the nerve of the contractible
groupoid with two objects $0$ and $1$.
\end{paragr}

\begin{definition}\label{def:Rezkmc}
We define the \emph{dendroidal Rezk model category} as the left
Bousfield localization of the model category for dendroidal
Segal spaces (\ref{defSegal}) by the maps
$$\Omega[T]\otimes J_d\To \Omega[T] \ , \quad T\in \Omega\, ,$$
obtained by tensoring with the unique morphism $J_d\To \eta$,
the image under $i_!$ of the unique map $J\to\Delta[0]$.
The fibrant objects of the dendroidal Rezk model category will be called
\emph{complete dendroidal Segal spaces}.
The weak equivalences of this model category structure will be called the
\emph{complete weak equivalences}.
\end{definition}

\begin{prop}\label{complinterval}
For any normal dendroidal set $A$, the map
$A\otimes J_d\To A$ induced by the projection $J_d\To \eta$ is a complete weak equivalence.
\end{prop}

\begin{proof}
It is sufficient to prove that this map is a weak equivalence
for the induced model category structure on $\sdset/A$, the latter being
equivalent to the category of presheaves on the
category $\Delta\times \Omega/A$ (where $\Omega/A$
is the category of elements of $A$). On the other hand, we know that
tensoring by $J_d$ preserves colimits as well as normal monomorphisms.
As, by virtue of \cite[Corollary 1.7]{dend3}, $\Omega/A$ is then a regular skeletal category
in the sense of \cite[8.2.3]{Ci3}, this proposition is a straightforward application of \cite[8.2.14]{Ci3}.
\end{proof}

\begin{cor}\label{cor:gentrivcofRezk000}
The inclusion $\dset\subset\sdset$ sends the weak operadic equivalences
between normal dendroidal sets to complete weak equivalences.
\end{cor}

\begin{proof}
By virtue of Propositions \ref{Segalspinesinnerhorns}
and \ref{complinterval}, this inclusion functor sends $J$-anodyne extensions in the
sense of \cite[3.2]{dend3} to complete weak equivalences.
This corollary thus immediately follows from \cite[Proposition 3.16]{dend3}
and from Ken Brown's Lemma \cite[Lemma 1.1.12]{Ho}.
\end{proof}

\begin{cor}\label{cor:gentrivcofRezk0}
Let $K$ be a generating set of trivial cofibrations in $\dset$.
We assume that all the maps in $K$ are morphisms between normal dendroidal sets
(which is a harmless hypothesis, by virtue of \cite[Remark 3.15]{dend3}).
A simplicial dendroidal set $X$ is a complete dendroidal Segal space
if and only if it is a Segal space, and if the map from $X$ to the terminal
simplicial dendroidal set has the right lifting property with respect to the inclusions
of shape
$$\partial\smp{n}\times B\cup \smp{n}\times A\To
\smp{n}\times B$$
for $j:A\To B$ in $K$ and $n\geq 0$.
\end{cor}

\begin{proof}
Let $L$ be the set of maps
$\partial\smp{n}\times B\cup \smp{n}\times A\To
\smp{n}\times B$, for $A\To B$ in $K$ and $n\geq 0$,
and consider the left Bousfield localization of the model
category for dendroidal Segal spaces by $L$.
It is clear that the $L$-fibrant objects are precisely
the Segal spaces $X$ such that the map from $X$
to the terminal object has the right lifting property
with respect to the elements of $L$.
Therefore, it is sufficient to prove that this left
Bousfield localization at $L$ coincides with the dendroidal
Rezk model category structure on $\sdset$.
Using corollaries \ref{simplnormalmono} and \ref{cor:gentrivcofRezk000},
we easily see that
the elements of $L$ are cofibrations and complete weak equivalences.
On the other hand,
for any trivial cofibration $X\to Y$ in $\dset$ and any
integer $n\geq 0$ the map
$\partial\smp{n}\times Y\cup \smp{n}\times X\To\smp{n}\times Y$
is a trivial cofibration of the localized model structure at $L$
(this readily follows from \cite[Lemma 4.2.4]{Ho}, for instance).
In particular, for any tree $T$, and $\varepsilon\in\{0,1\}$,
the map
$\Omega[T]=\Omega[T]\otimes\{\varepsilon\}\to\Omega[T]\otimes J_d$
is a weak equivalence of the localized model structure at $L$.
Therefore, any complete weak
equivalence is a weak equivalence for the left Bousfied localization by $L$.
\end{proof}

\begin{thm}\label{comploccstRezk}
The locally constant model category structure on $\sdset$ (\ref{def:loccstmc})
and the dendroidal Rezk model category structure (\ref{def:Rezkmc}) are equal.
\end{thm}

\begin{proof}
As these two model category structures on $\sdset$ have the same
class of cofibrations, it is sufficient to observe that they have the same class
of fibrant objects, which follows from
Corollary \ref{cor:gentrivcofRezk0} and from the characterization of fibrant objects
given by \ref{generators}, \ref{locconstReedyfibrant}~(i), \ref{dendroidalReedy}, and
\ref{chardendSegal}~(iii).
\end{proof}

\begin{cor}\label{inclQuillenequiv}
The inclusion functor $\dset\subset\sdset$ is a left Quillen
equivalence from the model category for $\infty$-operads
to the dendroidal Rezk model category.
\end{cor}

\begin{proof}
This follows from the preceding theorem and from
Proposition \ref{trivQuillenequivloccst}.
\end{proof}

\begin{cor}\label{cor:gentrivcofRezk00}
A morphism of dendroidal sets is a weak operadic equivalence if and only if its
image under the inclusion $\dset\subset\sdset$ is a complete weak equivalence.
\end{cor}

\begin{proof}
This is a reformulation the last assertion of Proposition \ref{trivQuillenequivloccst}
through Theorem \ref{comploccstRezk}.
\end{proof}

\begin{cor}
Let $X\To Y$ be a morphism between complete dendroidal Segal spaces. The following
conditions are equivalent.
\begin{itemize}
\item[(a)] The map $X\To Y$ is a complete weak equivalence.
\item[(b)] For any integer $n\geq 0$, the map $X_n\To Y_n$ is an equivalence of
$\infty$-operads.
\item[(c)] The map $X_0\To Y_0$ is an equivalence of $\infty$-operads.
\item[(d)] For any tree $T$, the map $X(T)\to Y(T)$ is a homotopy equivalence
between Kan complexes.
\end{itemize}
\end{cor}

\begin{paragr}
Consider the cosimplicial dendroidal set
\begin{equation}
\Delta_J:\Delta\To \dset
\end{equation}
defined by
\begin{equation}
\Delta_J[n]=i_!N\pi_1(\Delta[n])
\end{equation}
(so that $\Delta_J[1]=J_d$).
This cosimplicial object defines a unique colimit preserving functor
\begin{equation}\label{def:real}
\sdset\To \dset\ , \quad X\longmapsto \real X
\end{equation}
such that
\begin{equation}
\real{\Delta[n]\times\Omega[T]}=\Delta_J[n]\otimes\Omega[T]\, .
\end{equation}
The functor $\real{-}$ has a right adjoint
\begin{equation}\label{def:sing}
\dset\To\sdset\ , \quad X\longmapsto \sing X
\end{equation}
defined by
\begin{equation}
\sing{X}(T)_{n}=\Hom_\dset(\Delta_J[n]\otimes\Omega[T],X)\, .
\end{equation}
\end{paragr}

\begin{prop}\label{realleftQuillenequiv}
The functor \eqref{def:real} is a left Quillen equivalence from the
dendroidal Rezk model category to the model category for $\infty$-operads.
\end{prop}

\begin{proof}
Using the fact that $\dset$ is a monoidal model category,
it is easily seen that \eqref{def:real} is a left Quillen functor
from the generalized Reedy model structure (given by
Proposition \ref{dendroidalReedy}) to the model category for $\infty$-operads.
Therefore, to prove that this is a left Quillen functor
for the dendroidal Rezk model structure, it is sufficient to prove that
it sends inner horns as well as maps of shape 
$$\Omega[T]\otimes J_d\To \Omega[T] \ , \quad T\in \Omega\, ,$$
to weak operadic equivalences. But this latter property
follows from the fact that the composition of \eqref{def:real}
with the inclusion $\dset\subset\sdset$ is (isomorphic to) the
identity. Similarly, to prove that \eqref{def:real} is a left Quillen equivalence,
by virtue of Corollary \ref{inclQuillenequiv}, it is
sufficient to prove that its composition with
the inclusion $\dset\subset\sdset$ is a left Quillen
equivalence, which is more than obvious.
\end{proof}

\begin{paragr}\label{paragr:defK}
Let $\infty\text{-}\oper$ be the full subcategory of $\dset$
spanned by $\infty$-operads (\emph{i.e.} fibrant objects).
We define a functor
\begin{equation}\label{def:K}
K:\infty\text{-}\oper\To \sdset
\end{equation}
by the formula below (see \ref{reminderHoms}):
\begin{equation}
K(X)(T)_{n}=k(\Omega[T],X)_n\, .
\end{equation}
\end{paragr}

\begin{prop}\label{Knice}
The functor \eqref{def:K} takes its values in the full subcategory of
$\sdset$ spanned by complete dendroidal Segal spaces.
Moreover, it preserves fibrations as well as weak equivalences between
$\infty$-operads, and, under the canonical equivalence $\ho(\infty\text{-}\oper)\simeq\ho(\dset)$,
the corresponding functor
$$K:\ho(\dset)\To \ho(\sdset)$$
is canonically isomorphic to the functor
$$\derR\mathit{Sing}_J:\ho(\dset)\To \ho(\sdset)$$
(which is an equivalence of categories).
\end{prop}

\begin{proof}
In view of the identification of Proposition \ref{kAXMapAX},
this is a straightforward application of the general properties
of mapping spaces; see \cite[Propositions 5.4.1, 5.4.3 and 5.4.7]{Ho}
(remark that, for any normal dendroidal set $A$, $\Delta_J[\bullet]\otimes A$
provides a canonical cosimplicial frame of $A$ in the sense of
\cite[Definition 5.2.7]{Ho}).
\end{proof}

\begin{rem}\label{rem:BVtenrsorRezkcmc}
The Boardman-Vogt tensor product on $\dset$ induces a symmetric monoidal
structure on $\sdset$: for two simplicial dendroidal sets $X$ and $Y$,
their tensor product $X\otimes Y$ is simply defined termwise:
$$(X\otimes Y)_n=X_n\otimes Y_n\ , \quad n\geq 0\, .$$
Using the fact that $\dset$ is a symmetric monoidal model
category with respect to the Boardman-Vogt tensor product (Theorem \ref{cmfdendsetjoyalike}),
it is easily seen that $\sdset$, endowed with
the dendroidal Rezk model structure, is also a symmetric monoidal model category. Moreover,
the functor $\dset\subset\sdset$ is a symmetric monoidal left Quillen equivalence.
\end{rem}

\section{Segal pre-operads}

\begin{definition}
A \emph{Segal pre-operad} is a dendroidal space $A$ such that $A(\eta)$ is a discrete
simplicial set (\emph{i.e.} all the simplices of positive dimension in $A(\eta)$ are degenerated).
We denote by $\preoper$ the full subcategory of $\sdset$
spanned by Segal pre-operads.
\end{definition}

\begin{paragr}\label{existdiscretization}
The category of Segal pre-operads is in fact the category of presheaves on the category
$S(\Omega)$, which is obtained as the localization of $\Delta\times \Omega$
by the arrows of shape $([n],\eta)\To([m],\eta)$. We denote by
\begin{equation}
\gamma:\Delta\times \Omega\To S(\Omega)
\end{equation}
the localization functor. Under the identification $\preoper\simeq \pref{S(\Omega)}$, the
inverse image functor
\begin{equation}
\gamma^*:\preoper\To\sdset
\end{equation}
is simply the inclusion functor. The inclusion functor $\gamma^*$ thus has a
right adjoint
\begin{equation}
\gamma_*:\sdset\To\preoper
\end{equation}
as well as a left adjoint
\begin{equation}
\gamma_!:\sdset\To\preoper\, .
\end{equation}
The explicit description of these adjoints will be needed later on.

The right adjoint, $\gamma_*:\sdset\To\preoper$,
is defined as follows. Let $X$ be a dendroidal space. Then $\gamma_*(X)$
is the subobject of $X$ given by all the dendrices whose vertices
are degenerated. More explicitely, for a tree $T$, let write $E(T)$ for its
set of edges (colours), with the evident inclusion (which is natural in $T$)
$$v_T:\coprod_{e\in E(T)}\eta \To \Omega[T]\, .$$
For a simplicial set $K$, we shall identify the set $K_0$
with the corresponding discrete simplicial set, and write
$s:K_0\To K$ for the inclusion. Then $\gamma_*(X)(T)$ is defined
as the following pullback of simplicial sets.
\begin{equation}\label{descriptrightadj}\begin{split}
\xymatrix{
\gamma_*(X)(T)\ar[r]\ar[d]& X(T)\ar[d]^{v^*_T} \\
\prod_{e\in E(T)}X(\eta)_0\ar[r]^s&\prod_{e\in E(T)}X(\eta)
}\end{split}
\end{equation}

The left adjoint $\gamma_!:\sdset\To\preoper$ can also be made explicit as follows.
For a simplicial dendroidal set $X$ as above, consider the set $\pi_0X(\eta)$ of connected components
of the simplicial set $X(\eta)$. We have $\gamma_!(X)(T)=X(T)$ for any tree $T$ such that
there is no map $T\To\eta$ in $\Omega$. If there is a map $\varepsilon:T\To\eta$ in $\Omega$,
then its unique (remember there is a canonical isomorphism $\Omega/\eta= \Delta$),
and we can describe $\gamma_!(X)(T)$ as the pushout below.
\begin{equation}\label{descriptleftadj}\begin{split}
\xymatrix{
X(\eta)\ar[r]\ar[d]&X(T)\ar[d]\\
\pi_0X(\eta)\ar[r]&\gamma_!(X)(T)
}\end{split}
\end{equation}
\end{paragr}

\begin{paragr}
A morphism of Segal pre-operads is a monomorphism if and only if its image by $\gamma^*$ is
(because $\gamma^*$ is a fully faithful limit preserving functor).
We say that a morphism of Segal pre-operads $X\To Y$ is a \emph{normal monomorphism} if its image
by $\gamma^*$ has the same property (this just means that the map $X_n\To Y_n$
is a normal monomorphism of dendroidal sets for any integer $n\geq 0$).

A Segal pre-operad $X$ is \emph{normal} if $\varnothing\To X$ is a normal monomorphism.

A morphism of Segal pre-operads is a \emph{trivial fibration} if it has the right lifting property
with respect to the class of normal monomorphisms.
\end{paragr}

\begin{lemma}\label{lemma:leftadjnormmono}
If $X\To Y$ is a normal monomorphism of simplicial dendroidal sets and if
$\pi_0X(\eta)\To \pi_0Y(\eta)$ is injective, then $\gamma_!(X)\To\gamma_!(Y)$
is a normal monomorphism of Segal pre-operads.
\end{lemma}

\begin{proof}
One sees easily from the explicit description of $\gamma_!$ given by the pushouts
\eqref{descriptleftadj} that, for any tree $T$ above $\eta$, the map
$\gamma_!(X)(T)\To\gamma_!(Y)(T)$ is injective. For any tree $T$ which has a non
trivial automorphism in $\Omega$, there is no map from $T$ to $\eta$.
As, for such a tree $T$, we have $\gamma_!(X)(T)=X(T)$, it is clear that the map
$\gamma_!(X)\To\gamma_!(Y)$ is a normal monomorphism.
\end{proof}

\begin{prop}\label{prop:gennormmonoSeg}
Let $I$ be the set of maps
\begin{equation}\label{eq:gennormmonoSeg}
\gamma_!\big( \partial\smp{n} \times\Omega[T]\cup \smp{n}\times\partial\Omega[T]\big)
\To \gamma_!\big( \smp{n}\times\Omega[T]\big) 
\end{equation}
for any tree $T$ \emph{with at least one vertex}, and for any integer $n\geq 0$,
together with the map $\varnothing\To\eta$.
Then the smallest class of maps in $\preoper$ which is closed under
pushouts, transfinite compositions and retracts, and which contains $I$, is the
class of normal monomorphisms.
\end{prop}

\begin{proof}
Let us call \emph{$I$-cofibrations} the elements of the smallest class of maps
which contains $I$ and is closed under pushouts, transfinite compositions, and retracts.

If $T$ is a tree with at least one vertex, then, for any integer $n\geq 0$, the evaluation of the map
\begin{equation}\label{xxxxx}
 \partial\smp{n} \times\Omega[T]\cup \smp{n}\times\partial\Omega[T]\To\smp{n}\times\Omega[T] 
\end{equation}
at $\eta$ is bijective, so that, by virtue of Lemma \ref{lemma:leftadjnormmono},
its image by $\gamma_!$ is a normal monomorphism.
Hence any map in $I$ is a normal monomorphism of Segal pre-operads.
Therefore, any $I$-cofibration is a normal monomorphism.

Conversely, consider a normal monomorphism of Segal pre-operads $u:A\To B$.
Let $A'$ be the Segal pre-operad obtained from the pushout below.
$$\xymatrix{
\varnothing \ar[r]\ar[d]& A \ar[d]\\
\coprod_{b\in (B(\eta)_0-A(\eta)_0)} \eta \ar[r]& A'
}$$
Then the inclusion $A\To A'$ is certainly an $I$-cofibration, and one checks easily
that the canonical map $A'\To B$ is still a normal monomorphism.
Thus, to prove that $u:A\To B$ is an $I$-cofibration, we may assume,
without loss of generality, that the map $A(\eta)\To B(\eta)$ is bijective on $0$-simplices.
Applying the small object argument to the map $u$ with the set of maps \eqref{eq:gennormmonoSeg}
(for any tree $T$ with at least one vertex, and for any integer $n\geq 0$),
we obtain a factorization of $u$ of shape
$$A\xrightarrow{ \ v \ }C\xrightarrow{ \ p \ } B\, ,$$
in which $v$ is an $I$-cofibration, while $p$ has the right lifting property
with respect to maps of shape \eqref{eq:gennormmonoSeg}
(still for any tree $T$ with at least one vertex, and for any integer $n\geq 0$).
Moreover, one checks that $v$ induces a bijection by evaluating at $\eta$, which implies
that $p$ has the same property.
We claim that $\gamma^*(p)$ has the right lifting property with respect to maps of shape
\eqref{xxxxx} (for any tree $T$ and any integer $n$). Indeed, 
in the case $T$ has at least one vertex this follows
by a standard adjunction argument. In the case where $T=\eta$, this lifting property means that
the map $C(\eta)\To B(\eta)$ is a trivial fibration between discrete simplicial sets,
i.e. is a bijective map on the $0$-simplices. Hence, since $\gamma^*$ is
fully faithful, the map $p$ has the right lifting with respect to $u$.
By the retract argument \cite[Lemma 1.1.9]{Ho}, this implies that $u$ is a retract of $v$, whence is
an $I$-cofibration.
\end{proof}

\section{Segal operads}\label{segaloperads}

\begin{definition}\label{def:Segaloper}
A \emph{Segal operad} is a Segal pre-operad $X$ such that, for any tree $T$, the map
$$X(T)=X^{\Omega[T]}\To X^{\spine[T]}$$
is a trivial fibration of simplicial sets, where, if $A$ is a dendroidal set, $X^A$ denotes the
simplicial set whose $n$-simplices are the maps of dendroidal sets from $A$ to $X_n$
(with the notations of \ref{paragr:defexpA}, we thus have $\gamma^*(X)^A=X^A$).
We write $\segoper$ for the full subcategory of $\preoper$ spanned by
Segal operads.

A \emph{Reedy fibrant Segal operad} is a Segal pre-operad whose image by $\gamma^*$
is fibrant in the model category structure for dendroidal Segal spaces (see
Definition \ref{defSegal}). Note that any Reedy fibrant Segal operad is indeed a
Segal operad; see Corollary \ref{chardendSegal}.

A morphism of Segal pre-operads is a \emph{Segal weak equivalence} if its image
by $\gamma^*$ is a complete weak equivalence (\ref{def:Rezkmc}).

A morphism between Reedy fibrant Segal operads is \emph{fully faithful}
if its image by the functor $\gamma^*$ is fully faithful; see \ref{def:fullyfaithfulSegalspaces}.
\end{definition}

\begin{prop}\label{unitweakequiv}
Let $X$ be a dendroidal Segal space.
Then $\gamma_*(X)$ is a Reedy fibrant Segal operad.
\end{prop}

\begin{proof}
Since, for any tree $T$, the evaluation of the map $\spine[T]\To \Omega[T]$ at $\eta$
is bijective, the commutative square
$$\xymatrix{
\gamma^*\gamma_*(X)^{\Omega[T]}\ar[r]\ar[d]& X^{\Omega[T]}\ar[d]\\
\gamma^*\gamma_*(X)^{\spine[T]}\ar[r]& X^{\spine[T]}
}$$
is cartesian. Therefore, the functor $\gamma^* \gamma_*$ preserves
dendroidal Segal spaces. In other words, the functor $\gamma_*$ sends
dendroidal Segal spaces to Reedy fibrant Segal operads.
%
\end{proof}

\begin{lemma}\label{lemma:Segalnice0}
With the notations of paragraph \ref{paragr:defK},
for any $\infty$-operad $X$, the natural map $X\To K(X)$
is a complete weak equivalence.
\end{lemma}

\begin{proof}
As $X=K(X)_0$,
this is a reformulation of the fact that $K(X)$ is
a complete dendroidal Segal space; see
Proposition \ref{locconstReedyfibrant}, Theorem \ref{comploccstRezk}
and Proposition \ref{Knice}.
\end{proof}

\begin{prop}\label{prop:Segalnice}
For any dendroidal Segal space $X$, the map $X_0\To X$ is a
complete weak equivalence.
\end{prop}

\begin{proof}
Let $X$ be a dendroidal Segal space.
Given a bisimplicial object $U$, we write $\mathit{diag}(U)$ for the simplicial
object defined by $\mathit{diag}(U)_n=U_{n,n}$.
We define the bisimplicial dendroidal set $V$ by $V_{m,n}=X_n^{(\Delta[m])}$
(see paragraph \ref{reminderHoms}), and we put $W=\mathit{diag}(V)$.
The maps $\Delta[m]\To\Delta[0]$ induce embeddings $X_n=V_{0,n}\subset V_{m,n}$,
and thus a monomorphism $X\To W$.
Recall that $K(X_0)=V_{\bullet,0}$ is a fibrant resolution of $X_0$
in the dendroidal Rezk model structure; see Proposition \ref{Knice}.
We have a canonical commutative square of the following form
$$\xymatrix{
X_0\ar[r]\ar[d]&X\ar[d]\\
K(X_0)\ar[r]&W
}$$
in which the map $X_0\To K(X_0)$ is a complete weak equivalence (by Lemma \ref{lemma:Segalnice0}).
It is thus sufficient to prove that the maps $X\To W$ and $K(X_0)\To W$ are complete
weak equivalences. 

By virtue of Lemma \ref{lemma:Segalnice0},
the inclusion $X_n\to V_{\bullet,n}=K(X_n)$
is a weak equivalence for any integer $n\geq 0$.
Using Theorem \ref{comploccstRezk} (so that we can compute homotopy
colimits in $\sset^{\op\Omega}$ in the usual way), this implies that the induced map
$$X=\underset{\Delta[n]\in\op\Delta}{\mathrm{hocolim}}\, X_n
\To \underset{\Delta[n]\in\op\Delta}{\mathrm{hocolim}}\, K(X_n)=W$$
is a complete weak equivalence.

If we work with the projective model structure on $\sdset=\dset^{\op\Delta}$
associated to the model structure on $\dset$ (that is the model
category whose weak equivalences (or fibrations) are the maps
whose evaluation at each object of $\Delta$ is a weak
equivalence (or a fibration, respectively) in $\dset$), then,
for any tree $T$, the functor
$$\Delta[n]\mapsto\Delta[n]\times\Omega[T]$$
is a cosimplicial resolution of $\Omega[T]$, while
$$\Delta[m]\mapsto X^{(\Delta[m])}=V_{m,\bullet}$$
is a simplicial resolution of $X$.
Therefore, by virtue of \cite[Proposition 5.4.7]{Ho},
for any tree $T$, we can identify the simplicial set $W(T)$ with the
mapping space $\mathit{Map}(\Omega[T],X)$.
On the other hand, as the evaluation at zero functor $X\mapsto X_0$
is a right Quillen functor from $\dset^{\op\Delta}$ to $\dset$,
we have the following natural identifications in the
homotopy category of Kan complexes:
$$\mathit{Map}(\Omega[T],X_0)\simeq\mathit{Map}(\Omega[T],X)\, .$$
In other words, with the notations introduced in paragraph \ref{paragr:defK},
by virtue of Proposition \ref{kAXMapAX}, the map $K(X_0)\To W$
induces a canonical isomorphism
$$K(X_0)\simeq W$$
in the homotopy category of $\sset^{\op\Omega}$ (corresponding to termwise
weak equivalences of $\sset$). Therefore the map $K(X_0)\To W$
is a complete weak equivalence, and this ends the proof.
\end{proof}

\begin{rem}
If we keep the notations used in the proof above, we may use the
simplicial dendroidal set $W$ to obtain a canonical resolution of $X$ by
a complete dendroidal Segal space: one may consider a fibrant resolution
$W\To Y$ for the generalized Reedy model structure on $\sset^{\op\Omega}$.
Then the map $X\To Y$ is a complete weak equivalence, and $Y$ is a
complete dendroidal Segal space.
\end{rem}

\begin{cor}\label{cor:Segalnice}
The functor $\gamma_*$ sends complete weak equivalences between
dendroidal Segal spaces to Segal weak equivalences, and, for any
dendroidal Segal space X, the map $\gamma^*\gamma_*(X)\To X$ is
a complete weak equivalence.
\end{cor}

\begin{proof}
It is clearly sufficient to prove the last assertion, which follows from
the fact that, by virtue of Propositions \ref{unitweakequiv}
and \ref{prop:Segalnice},
for any dendroidal Segal space $X$, there exists a commutative square
$$\xymatrix{
\gamma^*\gamma_*(X)_0\ar[r]\ar@{=}[d]&\gamma^*\gamma_*(X)\ar[d]\\
X_0\ar[r]&X
}$$
in which three hence all maps are weak equivalences.
\end{proof}

\begin{cor}\label{cor:Segalnice}
A morphism between dendroidal Segal spaces $X\To Y$ is a complete weak equivalence if and
only if $\gamma_*(X)\To \gamma_*(Y)$ is a Segal weak equivalence of Segal operads.
\end{cor}

\begin{paragr}
Given a dendroidal Segal space $X$, there is a canonical operad $\mathit{ho}(X)$
associated to it, whose set of colours is $X(\eta)_0$, and whose sets of maps are given by
$\pi_0(X(x_1,\ldots,x_n;x))$ (the fact that this defines an operad
can be proved using the explicit description of the operad associated to
an $\infty$-operad (see \cite[Proposition 6.14]{dend3}),
Corollary \ref{cor:Segalnice} (to reduce to the case of a complete dendroidal Segal space),
Proposition \ref{kAXMapAX}, as well as the Quillen equivalence
of Proposition \ref{realleftQuillenequiv}; however, it is not difficult to
understand this construction in elementary terms).
\end{paragr}

\begin{definition}
A morphism of dendroidal Segal spaces $X\To Y$ is \emph{essentially surjective}
if the morphism of operads $\mathit{ho}(X)\To \mathit{ho}(Y)$ is essentially surjective.
\end{definition}

\begin{rem}\label{remark:equivDendSegalspacesvsSegaloperads}
A morphism of dendroidal Segal spaces $X\To Y$
is fully faithful (see Definition \ref{def:fullyfaithfulSegalspaces})
and essentially surjective if and ony if the induced morphism
$\gamma^*\gamma_*(X)\To\gamma^*\gamma_*(Y)$
has the same property.
\end{rem}

\begin{thm}\label{thm:charequivsegalspaces}
Let $f:X\To Y$ be a morphism of dendroidal Segal spaces.
The following conditions are equivalent.
\begin{itemize}
\item[(a)] The map $f$ is a complete weak equivalence;
\item[(b)] The map $\gamma_*(f):\gamma_*(X)\To\gamma_*(Y)$ is a weak
equivalence of Segal operads.
\item[(c)] The map $\gamma^* \gamma_*(f):\gamma^* \gamma_*(X)\To\gamma^* \gamma_*(Y)$ is
fully faithful and essentially surjective.
\item[(d)] The map $f$ is fully faithful and essentially surjective.
\end{itemize}
\end{thm}

\begin{proof}
The equivalence between (a) and (b) follows from
Corollary \ref{cor:Segalnice}, while the equivalence between (c) and (d) is a tautology.
To prove the remaining equivalences, we will use Theorem
\ref{charactwoperequiv2} to deduce the equivalence between
conditions (a) and (d). Indeed, we may assume that $X$ and $Y$ are complete (by
Corollary \ref{cor:Segalnice}),
so that the evaluated maps $X(T)\To Y(T)$ may be identified with the maps of mapping spaces
$\mathit{Map}(\Omega[T],X)\To \mathit{Map}(\Omega[T],Y)$;
up to the Quillen equivalence $\dset\subset\sdset$ (see Corollary \ref{inclQuillenequiv}),
and using Proposition \ref{kAXMapAX}, condition (d) (resp. (a))
may be interpreted by saying that the map between the corresponding
$\infty$-operads is fully faithful and essentially surjective (resp. satisfies
condition (b) of Theorem \ref{charactwoperequiv1}). This completes the proof
of the theorem.
\end{proof}

\begin{lemma}\label{lemma:segaltrivfibrezktrivfib}
If a morphism of Segal pre-operads has the right lifting property with respect
to normal monomorphisms, then it is a Segal weak equivalence.
\end{lemma}

\begin{proof}
Consider first a normal resolution $E_\infty$ of the terminal dendroidal set
(i.e. a cofibrant resolution of the terminal dendroidal set for the model
category structure of Theorem \ref{cmfdendsetjoyalike}),
considered as a simplicially constant simplicial dendroidal set.
We may see $E_\infty$ as a Segal pre-operad, and it follows immediately from
Theorem \ref{comploccstRezk} that, for any Segal pre-operad $X$, the
projection $E_\infty\times X\To X$ is a Segal weak equivalence.
Moreover, as $E_\infty$ is even a normal Segal pre-operad, $E_\infty \times X$ is always
normal. Therefore, it is sufficient to prove that, if $p:X\To Y$ is a morphism of normal
Segal pre-operads which has the right lifting property with respect to normal monomorphisms of pre-operads,
then it is a Segal weak equivalence. Note that $J_d$ may be seen as Segal pre-operad, so that,
for any normal pre-operad $A$, $J_d\otimes A$ is still a pre-operad,
and, whenever $A$ is normal, the map
$$A\amalg A=(\{ 0 \}\amalg \{1 \})\otimes A\To J\otimes A$$
is a normal monomorphism; see Corollary \ref{simplnormalmono}
and \cite[Proposition 1.9]{dend3}. We may now finish the proof in the standard way: as $Y$ is
normal, the map $p$ admits a section $s:Y\To X$, and the commutative square
$$\xymatrix{
X\amalg X\ar[r]^(.6){(1_X,sp)}\ar[d]& X\ar[d]^p\\
J_d\otimes X\ar[r]_{p\pi}\ar@{..>}[ur]^h &Y
}$$
admits a lifting $h$ (where $\pi:J_d\otimes X\To X$ denotes the map induced by the map $J_d\to\eta$), which
means that the map $p$ is a $J_d$-homotopy equivalence, whence a weak equivalence.
\end{proof}

\begin{thm}\label{thm:modcatSegaloperads}
The category of Segal pre-operads is endowed with a left proper cofibrantly generated
model category structure whose weak equivalences are the Segal weak equivalences, and
whose cofibrations are the normal monomorphisms.
\end{thm}

\begin{proof}
The preceding lemma tells us
that any trivial fibration of Segal pre-operads is a Segal weak equivalence.
On the other hand, by virtue of Proposition \ref{prop:gennormmonoSeg},
the class of normal monomorphism is generated by a small set of maps.
The existence of this model category is thus a particular case of
J.~Smith's Theorem; see \cite[Theorem 1.7 and Proposition 1.18]{beke1}.
The left properness property follows immediately from its counterpart for the
Rezk model category structure.
\end{proof}

\begin{paragr}
The model category structure above will be called the \emph{Reedy-Segal model category
structure} on $\preoper$. We will always consider the category
of simplicial dendroidal sets as a model category with the Rezk model structure (\ref{def:Rezkmc}).
By construction, the functor $\gamma^*:\preoper\To\sdset$ is a left Quillen functor.
Our purpose is to prove that it is in fact a left Quillen equivalence, and that
the fibrant pre-operads are precisely the Reedy fibrant Segal operads.
\end{paragr}

\begin{thm}\label{QuillenequivSegal}
The functor $\gamma^*:\preoper\To\sdset$ is a left Quillen equivalence
from the model category for Segal operads to the model category for
complete dendroidal Segal spaces.
\end{thm}

\begin{proof}
As $\gamma^*$ is a fully faithful left Quillen functor which preserves and detects weak equivalences,
this follows immediately from Corollary \ref{cor:Segalnice}.
\end{proof}

\begin{rem}
Note that Segal pre-operads are closed under tensor product (as defined in
Remark \ref{rem:BVtenrsorRezkcmc}), and that the model category of Theorem
\ref{thm:modcatSegaloperads} is symmetric monoidal, in such a way that the left Quillen
functor of Theorem \ref{QuillenequivSegal} is symmetric  monoidal as well
(this is immediate from Remark \ref{rem:BVtenrsorRezkcmc}).
\end{rem}

\begin{thm}\label{charReedyfibrantsegaloper}
Let $X$ be a Segal pre-operad. The following conditions are equivalent:
\begin{itemize}
\item[(a)] $X$ is fibrant in the model structure of Theorem \ref{thm:modcatSegaloperads};
\item[(b)] $X$ is a Reedy fibrant Segal operad;
\item[(c)] $X$ is a retract of $\gamma_*(Y)$ for some complete dendroidal Segal space $Y$;
\item[(d)] $X$ is a retract of $\gamma_*(Y)$ for some dendroidal Segal space $Y$.
\end{itemize}
\end{thm}

\begin{proof}
Condition (a) implies condition (b) because the inclusions
$\spine[T]\To \Omega[T]$ are trivial cofibrations in the symmetric monoidal
model category structure of Theorem \ref{thm:modcatSegaloperads}.

Let us prove that condition (b) implies condition (c).
If $X$ is Reedy fibrant, then we can choose a trivial cofibration
$\gamma^*(X)\To Y$ with $Y$ a complete dendroidal Segal space.
By virtue of Corollary \ref{cor:Segalnice}, we may assume
that the map $\gamma^*X\simeq \gamma^*\gamma_*\gamma^*(X) \To \gamma^*\gamma_*(Y)$
is a weak equivalence between fibrant objects in the model category for dendroidal Segal
spaces. As $\gamma^*$ is fully faithful, to prove that $X$ is a retract of $\gamma_*(Y)$,
it is sufficient to prove that the map $\gamma^*X \To \gamma^*\gamma_*(Y)$
is a cofibration (i.e. a normal monomorphism):
this follows from the fact that, by assumption, for any tree $T$, the group $\mathrm{Aut}(T)$
of automorphisms of $T$ in $\Omega$ acts freely on $Y(T)-X(T)$, and that we have
an $\mathrm{Aut}(T)$-equivariant inclusion of $\gamma_*(Y)(T)-X(T)$ in $Y(T)-X(T)$;
see the cartesian square \eqref{descriptrightadj}.

Condition (c) implies condition (a): as $\gamma_*$ is a right Quillen functor
(Corollary \ref{QuillenequivSegal}),
$\gamma_*(Y)$ is fibrant for any complete dendroidal Segal space $Y$, and the
class of fibrant objects of any model category is closed under retracts.

It is clear that condition (c) implies condition (d).

Finally, the fact that condition (d) implies condition (b)
follows from the fact that $\gamma_*$ sends dendroidal Segal spaces to
Reedy fibrant Segal operads (see the first assertion of Proposition \ref{unitweakequiv}).
\end{proof}

\begin{rem}
Recall from \ref{paragr:defK} the canonical functor
$$K:\infty\text{-}\oper\To \sdset\, .$$
We know that $K$ sends $\infty$-operads to complete dendroidal Segal spaces,
so that, by virtue of Proposition \ref{unitweakequiv}, we obtain a functor
$$\gamma_*K:\infty\text{-}\oper\To \segoper\, .$$
We also know from Proposition \ref{Knice} and from Theorem \ref{QuillenequivSegal}
that $\gamma_*K$ sends weak equivalences of $\infty$-operads to weak Segal equivalences, and that
the induced functor
$$\gamma_*K:\ho(\infty\text{-}\oper)\To\ho(\segoper)$$ is an equivalence
of categories.

Remark as well that any dendroidal set is a pre-operad, so that the inclusion
$\dset\subset\sdset$ factors through an inclusion $\dset\subset\preoper$
which happens to be a left Quillen equivalence (this follows immediately from
Corollary \ref{inclQuillenequiv} and from Theorem \ref{QuillenequivSegal}).
If $X$ is an $\infty$-operad, seen as a Segal pre-operad, then $\gamma_*K(X)$ is
a canonical fibrant replacement of $X$ in the model category of Theorem \ref{thm:modcatSegaloperads}.
\end{rem}

\begin{rem}\label{rem:references}
The identification $\dset/\eta=\sset$ allows us to deduce the Joyal model
category structure for quasi-categories from the model category structure for
$\infty$-operads; see \cite[Corollary 2.10]{dend3}.
Similarly, the dendroidal Rezk model structure of Definition \ref{def:Rezkmc}
induces Rezk's original model structure for complete Segal spaces (\cite[Section 12]{Rezk}),
while the model category structure for Segal categories can be obtained
from the model category structure of Theorem \ref{thm:modcatSegaloperads}
by slicing over $\eta$ as well. The Quillen equivalences
relating the homotopy theories of Segal categories and of
complete Segal spaces, proved by Joyal and Tierney
in \cite{joytier4} are deduced immediately from
their dendroidal analog, namely Corollary \ref{inclQuillenequiv}
and Proposition \ref{realleftQuillenequiv}, while the
Quillen equivalence from the model category for Segal categories
to the model category for complete Segal spaces, proved by
Bergner in \cite{bergner}, is a direct consequence of Theorem \ref{QuillenequivSegal}.
\end{rem}

\providecommand{\bysame}{\leavevmode\hbox to3em{\hrulefill}\thinspace}

\end{document}